\documentclass[10pt,a4paper]{article}
\usepackage[utf8]{inputenc}
\usepackage[english]{babel}
\usepackage{amsmath}
\usepackage{amsthm}
\usepackage{amsfonts}
\usepackage{amssymb}
\usepackage{graphicx}
\usepackage{a4wide}
\usepackage{bbm}
\usepackage{tikz}
\usepackage{todonotes}
\usetikzlibrary{calc}
\usepackage{cite}
\usepackage{authblk}
\usepackage{float}
\usepackage{subcaption}
\usepackage{filecontents}

\newcommand{\R}{\mathbb R}
\newcommand{\C}{\mathbb C}
\newcommand{\N}{\mathbb N}
\newcommand{\Z}{\mathbb Z}
\newcommand{\E}{\mathbb E}
\newcommand{\Prob}{\mathbb P}
\newcommand{\Cov}{\mathrm{Cov}}

\newcommand{\iid}{\stackrel{\text{i.i.d.}}{\sim}}
\newcommand{\X}{\mathcal X}
\newcommand{\Y}{\mathcal Y}

\newcommand{\comment}[1]{}

\DeclareMathOperator{\supp}{supp}

\DeclareMathOperator{\lspan}{span}
\DeclareMathOperator{\id}{id}
\DeclareMathOperator{\argmin}{arg\,min}

\DeclareMathOperator{\ran}{ran}
\DeclareMathOperator{\dom}{dom}

\DeclareMathOperator{\vol}{vol}

\newtheorem{theorem}{Theorem}[section]
\newtheorem{lemma}[theorem]{Lemma}
\newtheorem{assum}[theorem]{Assumption}

\newtheorem{corol}[theorem]{Corollary}

\numberwithin{equation}{section}

\newcommand{\widerhat}[1]{\mathpalette\widerhataux{#1}}
\newcommand{\widerhataux}[2]{\vphantom{#1\widehat{#2}}
	\tikz[baseline]{
		\node[inner sep=0, anchor=base] (widerhat){$#1#2$};
		\draw (widerhat.north west) ++ (.05em,.4ex) -- ($(widerhat.north west)!.5!(widerhat.north east)+(.05em,.8ex)$) -- ($(widerhat.north east)+(.05em,.4ex)$);
		}
	}

\usepackage{hyperref}
\hypersetup{
	colorlinks=true,
	urlcolor=green,
	linkcolor=red,
	citecolor=blue
	}

\title{Minimax detection of localized signals in statistical inverse problems}
\date{}

\author[1]{Markus Pohlmann}
\author[2]{Frank Werner}
\author[1,3,4,5]{Axel Munk}

\affil[1]{Institute for Mathematical Stochastics, University of Göttingen, Goldschmidtstr. 7, 37077 Göttingen, Germany}
\affil[2]{Institute for Mathematics, University of Würzburg, Emil-Fischer-Str. 30, 97074 Würzburg, Germany}
\affil[3]{Felix-Bernstein-Institute for Mathematical Statistics in the Biosciences, University of Göttingen, Goldschmidtstr. 7, 37077 Göttingen, Germany}
\affil[4]{Max Planck Institute for Biophysical Chemistry, Am Faßberg 11, 37077 Göttingen, Germany}
\affil[5]{Cluster of Excellence MBExC: ``Multiscale Bioimaging: From Molecular Machines to Networks of Excitable Cells", University Medical Center, Robert-Koch-Str. 40, 37075 Göttingen, Germany}

\begin{document}
\maketitle

\begin{abstract}
We investigate minimax testing for detecting local signals or linear combinations of such signals when only indirect data is available. Naturally, in the presence of noise, signals that are too small cannot be reliably detected. In a Gaussian  white noise model, we discuss upper and lower bounds for the minimal size of the signal such that testing with small error probabilities is possible. In certain situations we are able to characterize the asymptotic minimax detection boundary. Our results are applied to inverse problems such as numerical differentiation, deconvolution and the inversion of the Radon transform.
\end{abstract}

\textit{Keywords:} Hypothesis testing, minimax signal detection, statistical inverse problems; wavelet-vaguelette decomposition. 

\textit{AMS classification numbers: } 62F03, 65J22, 65T60, 60G15. \\

\section{Introduction}
\begin{figure}[ht]
	\centering
	\begin{subfigure}{.3\textwidth}
		\centering
		\includegraphics[width=\textwidth]{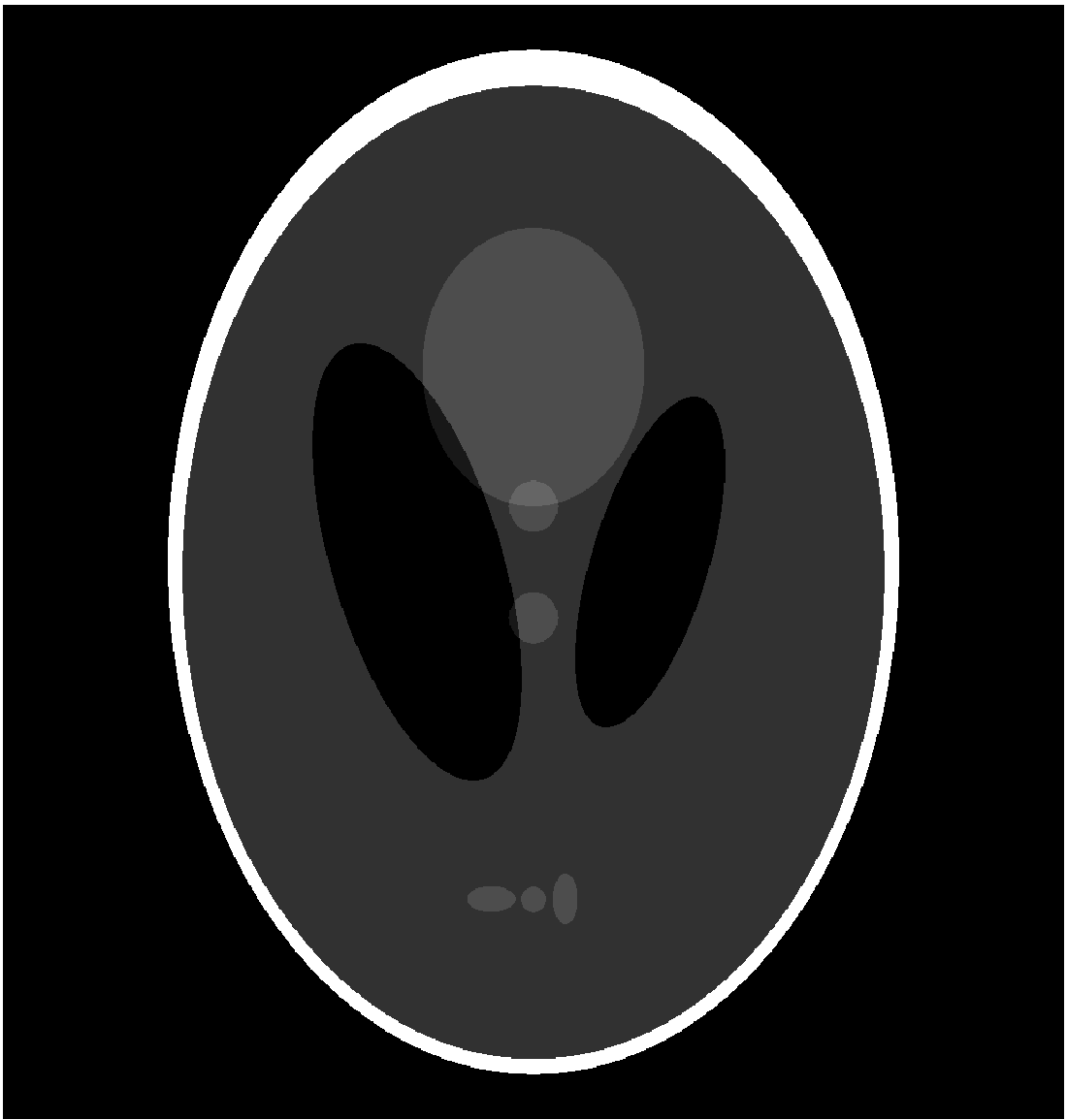}
		\caption{Reference image\label{subfig:refsignal}}
	\end{subfigure}
	\begin{subfigure}{.3\textwidth}
		\centering
		\includegraphics[width=\textwidth]{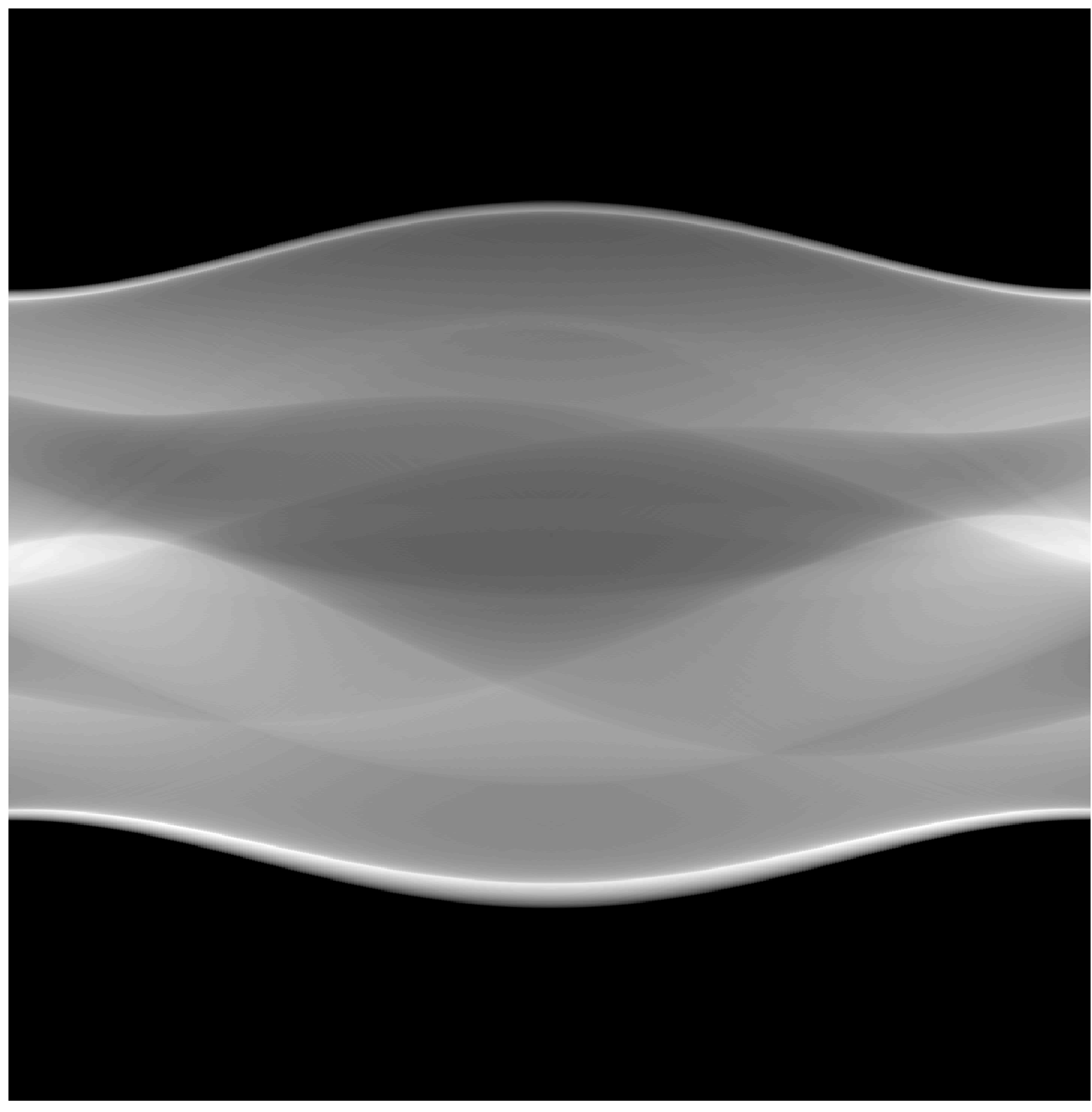}
		\caption{Sinogram of (a)\label{subfig:refsignalradon}}
	\end{subfigure}
	\begin{subfigure}{.3\textwidth}
		\centering
		\includegraphics[width=\textwidth]{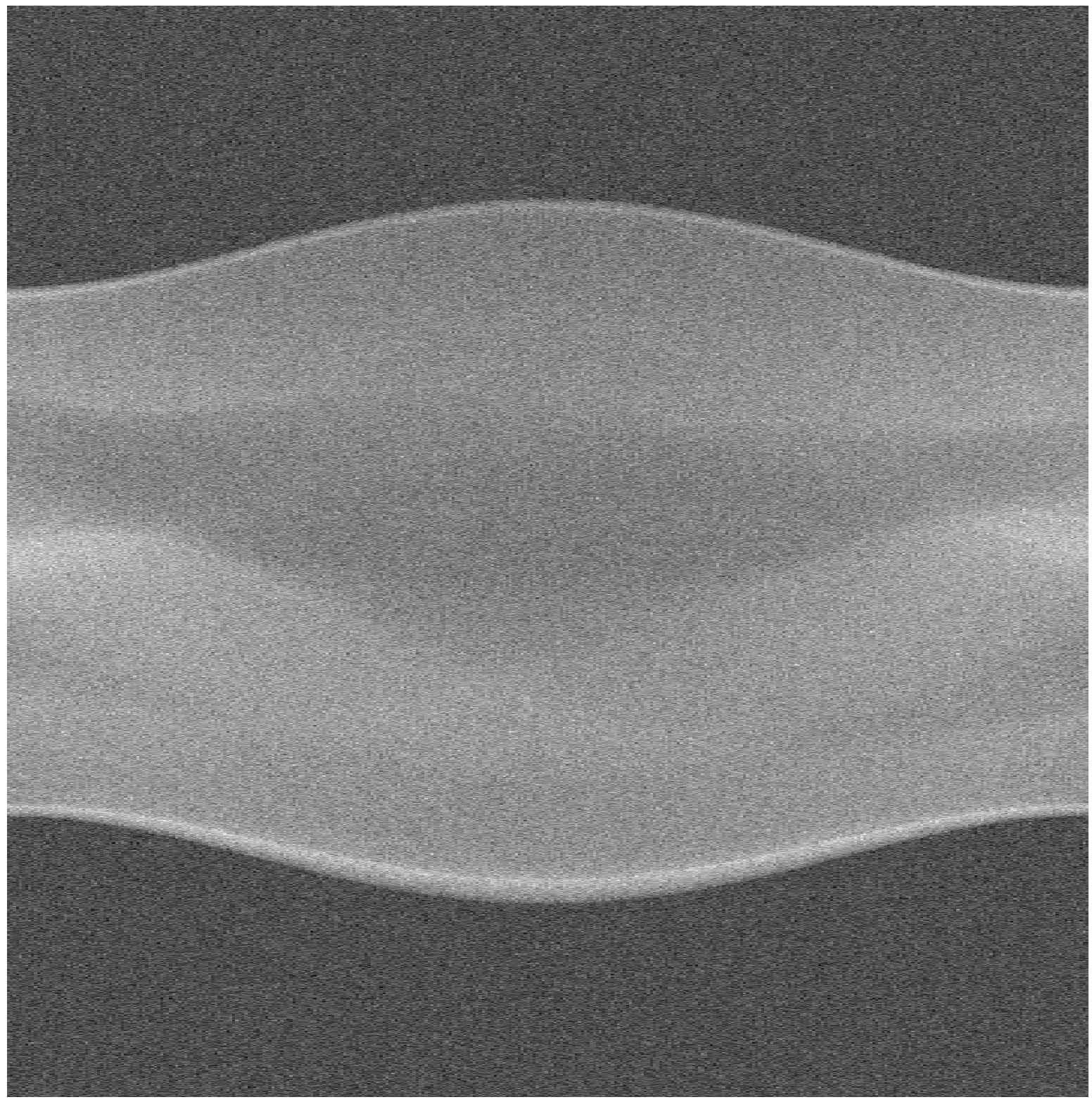}
		\caption{Noisy sinogram of (a)\label{subfig:refsignalradonerr}}
	\end{subfigure}\\
	\begin{subfigure}{.3\textwidth}
		\centering
		\includegraphics[width=\textwidth]{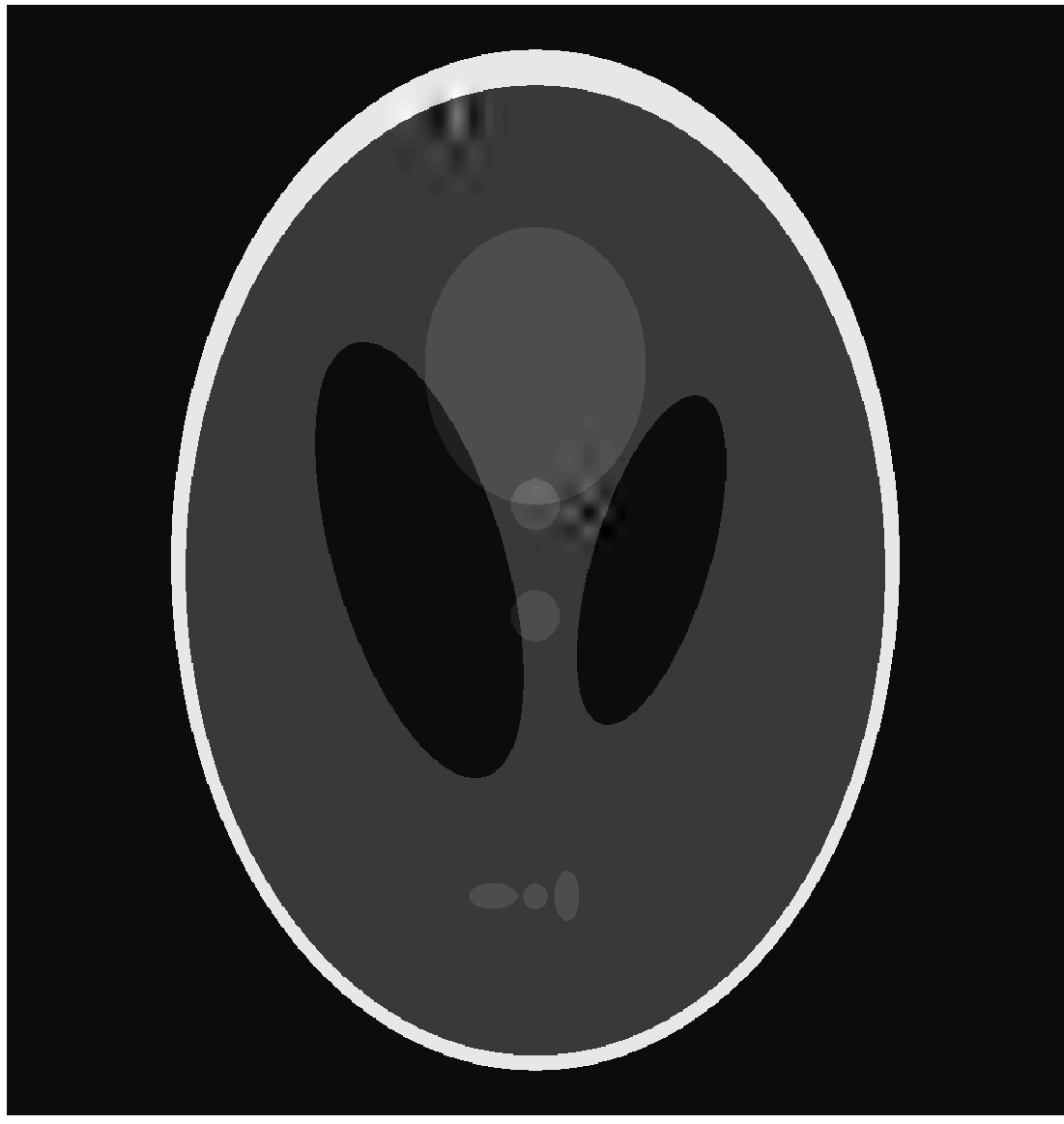}
		\caption{Distorted Image \label{subfig:refsignalanom}}
	\end{subfigure}
	\begin{subfigure}{.3\textwidth}
		\centering
		\includegraphics[width=\textwidth]{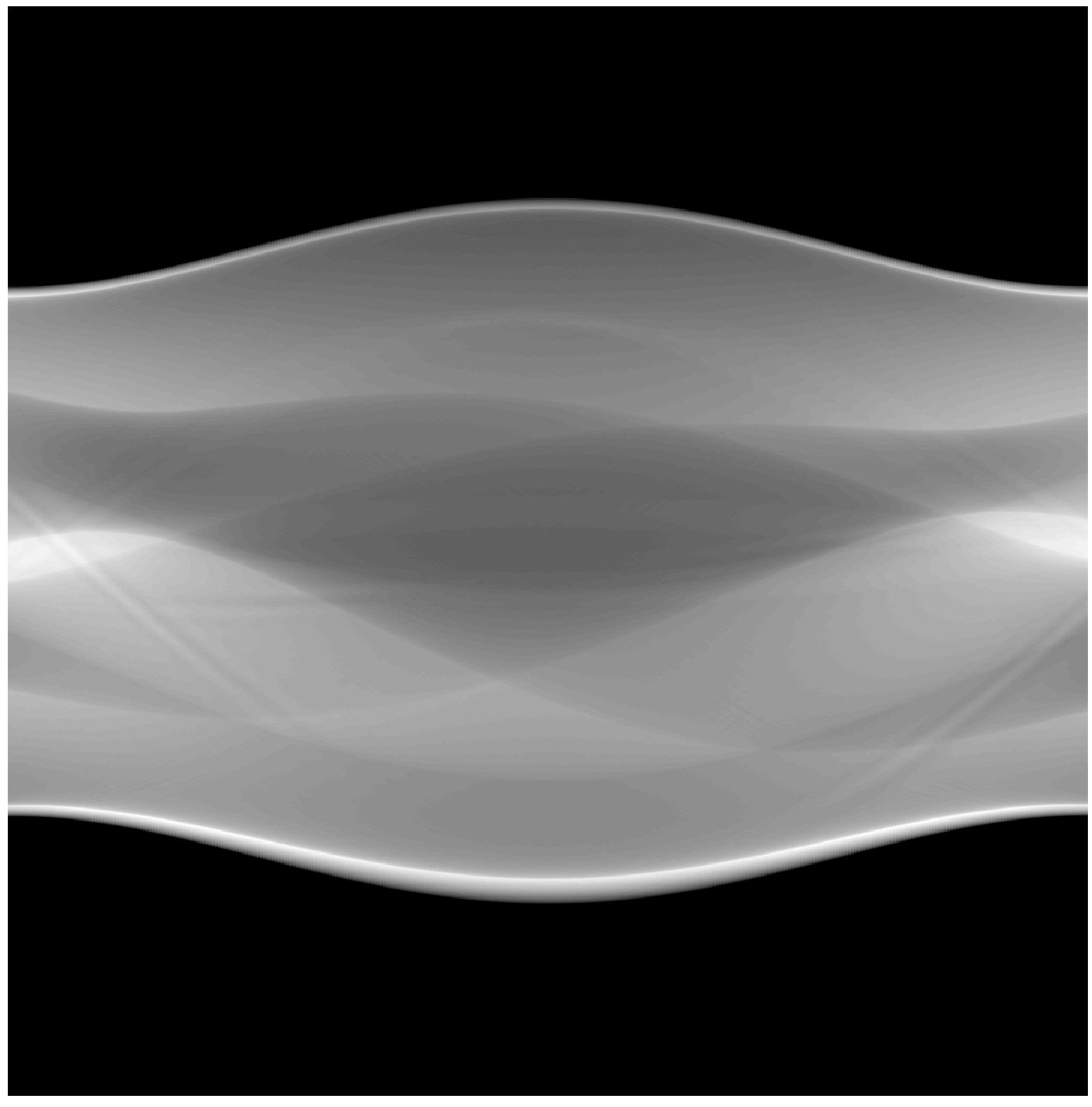}
		\caption{Sinogram of (d)\label{subfig:refsignalanomradon}}
	\end{subfigure}
	\begin{subfigure}{.3\textwidth}
		\centering
		\includegraphics[width=\textwidth]{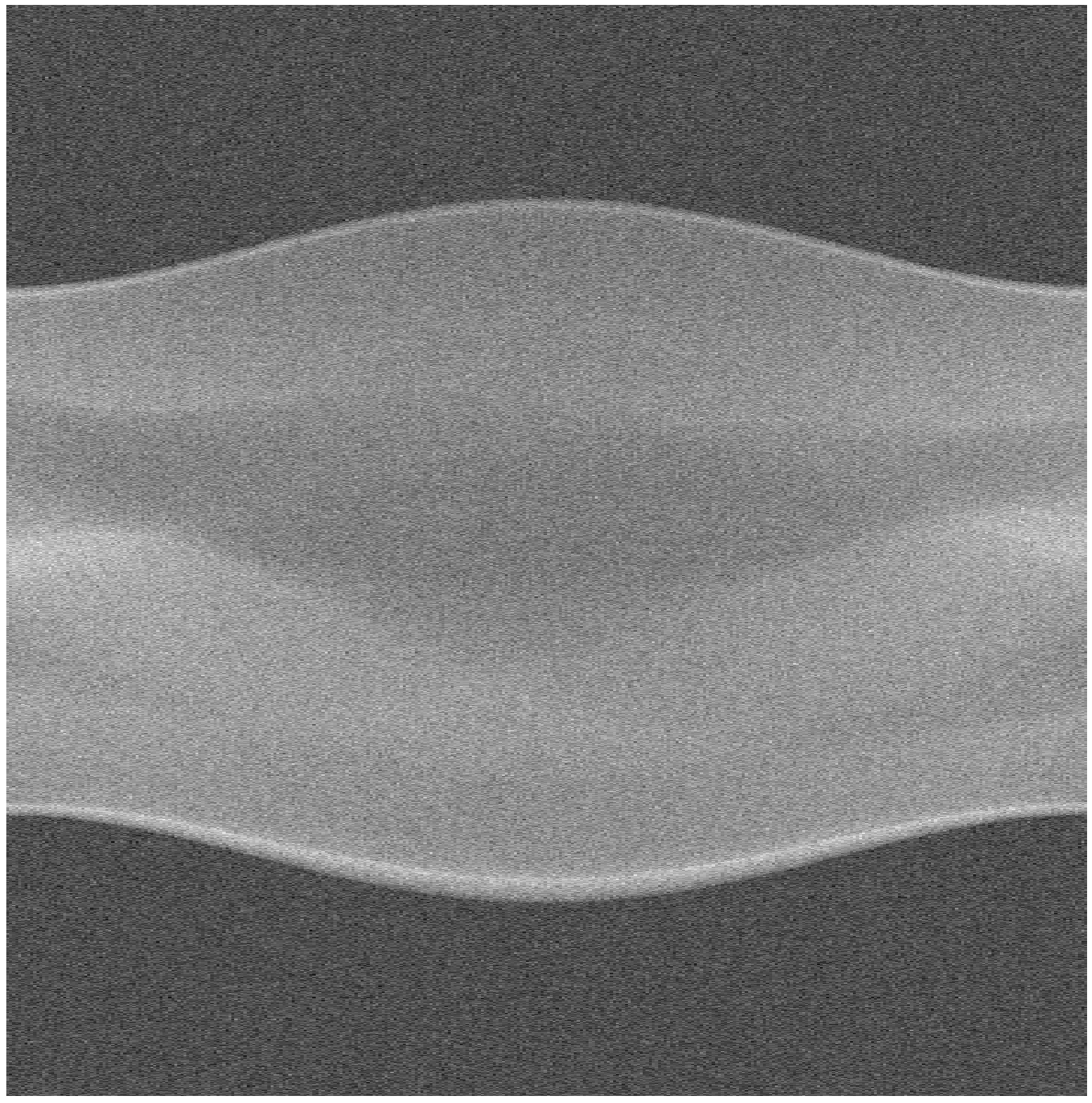}
		\caption{Noisy sinogram of (d)\label{subfig:refsignalanomradonerr}}
	\end{subfigure}
	\caption{Illustration of structured hypothesis testing in the CT example. To infer whether the unknown signal deviates from a reference image, we use a test based on the noisy sinogram. In the above example, when the distortion is assumed to be a linear combination of certain wavelets (cf. Sections \ref{subsec:case2_examples} and \ref{sec:simulations}), then the results of Theorem \ref{thm:case2_upperbound_nonasymp} imply the existence of a test which is able to distinguish the distorted (\ref{subfig:refsignalanom}) from the undistorted image (\ref{subfig:refsignal}) with type I and type II error both at most $0.05$, based on the measurements \ref{subfig:refsignalanomradonerr}.}\label{fig:radonintro}
\end{figure}

In many practical applications one aims to infer on properties of a quantity which is not directly observable. As a guiding example, consider computerized tomography (CT), where the interior (more precisely the tissue density) of the human body is imaged via the absorption of X-rays along straight lines. Mathematically, the relation between the available measurements $Y$ (absorption along lines, the so-called sinogram) and the unknown quantity of interest $f$ (the tissue density) is described by the Radon transform, which is an integral operator to be described in more detail later (cf. Figure \ref{fig:radonintro} for illustration). Potential further applications include astronomical image processing, magnetic resonance imaging, non-destructive testing and super-resolution microscopy, to mention a few. Typically, the measurements are either of random nature themselves (as e.g. in positron emission tomography (PET, see \cite{vardi1985}), magnetic resonance imaging (MRI, see \cite{klosowski2017}) or super-resolution microscopy (see \cite{munkstaudtwerner2020})) and/or additionally corrupted by measurement noise. This motivates us to consider the inverse Gaussian white noise model 
\begin{equation}\label{eq:inverseproblem_schematic}
	Y_\sigma= Af+\sigma \xi
\end{equation}
with a (known) bounded linear operator $A : \X \to \Y$ mapping between (real or complex) Hilbert spaces $\X$ and $\Y$, noise level $\sigma> 0$ and a Gaussian white noise $\xi$ on $\Y$ (details will be given in section \ref{sec:prelim}).

A major effort of research is devoted to the development and analysis of estimation and recovery methods of the signal $f$ from the measurements $Y_\sigma$ (see Section \ref{subsec:literature} for some references). However, when $f$ is expected to be very close to some reference $f_0$, by which we mean that either $f=f_0$ or $f$ deviates from $f_0$ by only a few localized components (anomalies), then instead of full recovery of $f$, one might be more interested in testing whether $f=f_0$ or not. This is especially relevant, since, when the signal-to-noise level is too small for full recovery, then testing may still be informative as it is well-known to be a simpler task (see e.g. \cite{proksch2018} and the references therein). Although of practical importance, testing in model \eqref{eq:inverseproblem_schematic} is a much less investigated endeavor than estimation and a full theoretical understanding has not been achieved yet. Hence, in this paper, we are interested in analyzing such testing methodology for inferring on $f$ based on the available data $Y_\sigma$. Note that, due to the linearity of the model \eqref{eq:inverseproblem_schematic}, we can w.l.o.g. assume that $f_0=0$. Thus, we suppose that either $f = 0$ (no anomaly is present) or $f= \delta u$ (an anomaly given by $\delta u$ is present), where $u\in \mathcal{F}_\sigma$ for some (finite) class $\mathcal{F}_\sigma\subseteq \X$ of non-zero functions, that are suitably normalized, and the constant factor $\delta$ describes its orientation, and -- more importantly -- how ``large" or ``pronounced" the signal $f$ is. {To this end, w}e consider the family of  testing problems
\begin{equation}\label{eq:testproblem_A}
	H_0: f=0 \quad\text{against}\quad H_{1,\sigma}: f = \delta u \text{ for some } u\in \mathcal{F}_\sigma \text{ and } |\delta|>\mu_\sigma,
\end{equation}
where $(\mu_\sigma)_{\sigma>0}$ is a family of {decreasing} non-negative real numbers. This can be viewed as the problem of \textit{detecting} an anomaly from the set $\{\delta u : u\in \mathcal{F}_\sigma, |\delta|>\mu_\sigma\}$. {Note that the sets $\mathcal F_\sigma$ of possible anomalies will become larger with smaller values of $\sigma$.}

{We would like to emphasize already at this point, that the type of anomalies we have in mind are \textit{local} deviations from zero, and hence, the influence of $A$ in \eqref{eq:inverseproblem_schematic} is expected to propagate to the testing problem \eqref{eq:testproblem_A} in an ill-posed way.}

We suppose that the family of classes $(\mathcal{F}_\sigma)_{\sigma>0}$ is chosen in advance. This choice is crucial for the analysis of the problem and it depends solely on the specific application: For CT we might think of small inclusions such as tumors, cf. Figure \ref{fig:radonintro}, where certain wavelets are used as mathematical representation. If no a priori knowledge about potential anomalies is known, it is natural to start by considering dictionaries $(u_k)_{k \in I}$ with good expressibility in $\X$, e.g. frames or wavelets, and set $\mathcal{F}_\sigma= \{u_k: k\in I_\sigma\}$ for subsets $I_\sigma$ of $I$. The particular choices that we analyze in this paper will be built from such dictionaries, see also \cite{ebner2020} and \cite{hubmer2021} for recent references in the context of estimation.

\subsection{Aim of the paper}

Given a family of classes $(\mathcal{F}_\sigma)_{\sigma>0}$, our main objective will be to assess to what extent powerful tests for the testing problem \eqref{eq:testproblem_A} exist. The answer will usually depend on the size of $\mu_\sigma$: If $\mu_\sigma$ is large enough, then powerful tests exist, and if $\mu_\sigma$ is too small, then no test has high power. Hence, we aim to find a minimal family of thresholds $(\mu_\sigma^*)_{\sigma>0}$, such that powerful detection at a controlled error rate is still possible. Vice versa, such a minimal family would determine which signals can not be detected reliably, even when they are present. 

To this end, we extend the existing theory on minimax signals detection in inverse problems focusing on localized signals and linear combinations of localized signals, which are common in practice. This has, to the best of our knowledge, not been investigated yet. We present upper bounds, lower bounds and asymptotics for the minimal values of $\mu_\sigma$ such that powerful tests for testing problems given by \eqref{eq:testproblem_A} exist. They depend on the difficulty of the inverse problem induced by the forward operator $A$, the cardinality of $I_\sigma$ (denoted by $|I_\sigma|$) and the inner products between the images $Au_k$, $k\in I_\sigma$, of the potential anomalies. {We aim to determine not just the asymptotic rate of $\mu_\sigma^*$, but also the corresponding minimax constant. Let us stress} that our results can be applied to a variety of dictionaries $(u_k)_{k \in I}$, such as wavelets, whereas previous results were restricted to dictionaries based on the SVD of the operator $A$. {This is a severe limitation as the shape of the signal in any real world application is fully unrelated to the operator (measurement device).} As one particular example, our results can be applied to the situation where the dictionary $(u_k)_{k \in I}$ is (a subset of) the famous Wavelet-Vaguelette-decomposition (WVD, see \cite {donoho1995wvd}) or the Vaguelette-Wavelet-decomposition (VWD, see \cite{abramovich1998vwd}) of $A$.

Figure \ref{fig:radonintro} serves as an illustrative example. If it is known a priori, that the anomaly which distorts the reference image is a linear combination of a certain collection of wavelets (see the discussion in Sections \ref{subsec:case2_examples} and \ref{sec:simulations} for details), then our results suggest that the anomaly that is present {in display (d)} is large enough, such that there is a test which is able to distinguish {it} from the undistorted image {in display (a)} with type I and type II error both at most $0.05$, based on the measurements in {display (f)} (see Theorem \ref{thm:case2_upperbound_nonasymp}){, despite the fact that the distortion is invisible by eye in display (f) compared to display (c)}. {In fact, the distrotion shown in Figure \ref{fig:radonintro} is deep in the alternative in the sense that Theorem \ref{thm:case2_upperbound_nonasymp} ensures its detectability by a test with high probability}. Note that our results are not restricted to wavelets. In fact, most of our results are applicable under very mild conditions on the dictionary $(u_k)_{k \in I}$.

We finnally stress that this paper does not constitute an exhaustive study of the subject. Rather, we aim to provide some first analysis and discuss some illustrative examples.

\subsection{Connection to existing literature}\label{subsec:literature}

{Most of the literature about inverse problems of the form \eqref{eq:inverseproblem_schematic} is concerned with the estimation of $f$ from $Y_\sigma$ (or its discretized version). 
	We mention the seminal monograph \cite{ehn96}, which utilizes for estimation of $f$ a spectral decomposition of $A$ in a deterministic noise model, and the more recent book \cite{hanke2017}.
	In case of random data, spectral estimation has been extensively treated in \cite{bhmr07}, and in \cite{m09} for the related problem of density deconvolution. 
	Whereas spectral methods are particularly well-suited for signals which are sparse in the spectral domain, estimation techniques which provide spatially sparse solutions include Wavelet thresholding \cite{donoho1995wvd,js97,abramovich1998vwd,ab06}, localized TV regularization \cite{am20}, and other nonlinear estimation schemes \cite{wernerhohage2012}. For a recent survey, see \cite{hlm22}. 
	Statistical minimax optimality \cite{abramovich1998vwd,am20,donoho1995wvd,johnstone2004,rs16,t00} and adaptation \cite{j99,l15} to the unknown smoothness of the signal are meanwhile well established theories, which resemble many of these methods as highly efficient (possibly after certain practical modifications).
	More recently, also Bayesian estimators \cite{gn20,mnp21} have been made accessible to a rigorous statistical analysis which theoretically supports their great practical success, known for a long time.} 

{In contrast, the paper at hand focuses on (minimax) \textit{testing}, which is well known to differ substantially from estimation. In case of the direct problem, i.e. when $\X=\Y$ and $A$ is the identity, meanwhile a comprehensive theory of minimax testing has been developed, see e.g. the seminal monograph \cite{ingster2003} based on a series of papers by Y. Ingster \cite{ingster1993}. These works treat a variety of problems to test the hypothesis ``$f=0$" against alternatives of the form ``$f\in \mathcal{F}_\sigma \text{ and } \Vert f\Vert_\X\geq\mu_\sigma$", where $\mathcal{F}_\sigma$ is a certain class of functions. Typical examples for $\mathcal F$ in \cite{ingster2003} are Sobolev or Besov balls $\mathcal F_\sigma$, i.e. $\mathcal F_\sigma$ is defined in terms of (global) smoothness properties of the anomaly $f$. In the last decades, this approach has been extended to the case of inverse problems (i.e. $A$ is allowed to differ from the identity) by several authors, which we will describe in more detail now. Early references on the closely related inverse problem of testing in density deconvolution models are \cite{hbm05,bmp09}, where the authors provide the rate behavior of $\mu_\sigma^*$ in different scenarios such as polynomially smooth kernels (noise densities) and different smoothness classes of the density. The suggested tests are based on kernel density estimators resulting from regularized Fourier inversion. Closely related to density deconvolution are error in variable models, where tests have been suggested based on similar approaches, see e.g. \cite{hbm05,ot21}. The first result on general inverse problems is \cite{laurent2011directindirect}, where the authors consider generic forward operators $A$ and $\X$-ellipsoids of the form
	\[
	\mathcal E_{a,2}^{\X} \left(R\right) = \left\{f \in \X: \sum_{j=1}^\infty a_j^2 \langle f, \xi_j\rangle_{\X}\leq R^2 \right\},
	\]
	where $\left(s_j,\xi_j, \eta_j\right)_{j \in \N}$ is a singular value decomposition (SVD) of $A$. In this case, it is shown that the inverse testing problem \eqref{eq:testproblem_A} with $\mathcal F_\sigma = \mathcal E_{a,2}^{\X} \left(R\right)$ is - in terms of the asymptotic rate  of $\mu_\sigma^*$ ensuring distinguishability, see Section 2.3 below - equivalent to a direct testing problem
	\[
	H_0^{\mathrm{DP}}:  Af = 0 \qquad\text{vs.}\qquad H_1^{\mathrm{DP}}: Af \in \mathcal E_{s^{-1}a,2}^{\Y}(R), \Vert Af\Vert_{\Y} \geq \mu_\sigma',
	\]
	with the sequence $s^{-1}a = (s_j^{-1}a_j)_{j \in \N}$. However, as also shown in \cite{lm14}, minimax testing procedures for the inverse testing problem \eqref{eq:testproblem_A} are not automatically minimax for the corresponding direct problem. For a detailed analysis of the relation between direct and inverse testing under global smoothness assumptions, we refer to \cite{ingster2012, ingster2014, marteau2014general, autin2019maxiset}, including extensions e.g. to multidimensional problems or specific tests.} {We also mention \cite{ot21} for a recent work in that spirit.}

{In contrast to the above mentioned work, the testing problem \eqref{eq:testproblem_A} we have in mind and our aim substantially differ from the previously investigated scenarios in two major aspects. First of all, instead of testing against a global smoothness condition $f \in \mathcal E_{a,2}^{\X} \left(R\right)$ with sufficiently large norm, we aim to test for localized anomalies, which are contained in a very specific set of candidate functions. An illustration is given in Figure 1, where the reference image in display (a) is distorted by a localized function in display (d). This localized function has certain smoothness properties (as a Wavelet), but has much more structure than just its smoothness. Secondly, it is not representable as a finite combination of singular vectors of the radon transform $A$ (which are non-local), which reveals SVD based methods as not well-suited. Nevertheless, it is compatible with the operator in the sense of a WVD. Therefore, we aim to investigate localized testing instead of global testing. This implies that in general it is not possible to transform the problem into a direct testing problem. A more detailed motivation for using such localized alternatives in inverse problems, and a discussion of the resulting difficulties can e.g. be found in \cite{proksch2018,kretschmann2022}.}

{For such testing problems, we aim to determine the exact asymptotic detection boundary of the problem \eqref{eq:testproblem_A} in the sense that we are not only interested in the decay rate of $\mu_\sigma^*$, but also in the constant which describes the phase transition between distinguishable and undistinguishable anomalies (alternatives). This paper is also motivated from previous work on minimax bump testing in time series \cite{enikeeva2018,enikeeva2020}, which to some extent can be rewritten as a one-dimensional inverse problem as in \eqref{eq:inverseproblem_schematic}, but to the best of our knowledge, no such result (neither for bump functions nor for other anomalies) is known in case of general inverse problems. We stress, that our work also shows that localized functions such as Wavelets turn out to be more suitable than bumps in inverse (testing) problems from a modeling and analysts view point, see e.g. \cite{donoho1995wvd} and our results below.}

{Finally, we want to highlight \cite{laurentloubes2012} explicitly, as this work considers alternatives consisting of linear combinations of anomalies given in terms of the SVD of the operator $A$, which is a special case of this study.}

\subsection{Outline}

We start by giving a detailed description about our model and some basic facts about testing and minimax signal detection in section \ref{sec:prelim}. Section \ref{sec:results} contains the main results: In section \ref{subsec:results_case1} we assume that $\mathcal{F}_\sigma$ is a collection of frame elements, and in section \ref{subsec:results_case2} we assume that $\mathcal{F}_\sigma$ contains functions in the linear span of a collection of frame elements. Both sections also include discussions about conditions that frames need to satisfy for our results to be applicable. We present illustrative simulation studies in section \ref{sec:simulations}. All proofs are postponed to section \ref{sec:proofs}.

\section{Preliminaries}\label{sec:prelim}

{Throughout the paper we assume, that the Hilbert space $\Y$ is separable, i.e. has a countable complete orthonormal system $\{ e_i : i \in \N\}$. This means, that
	\[
	\langle e_i, e_j\rangle_{\Y} = \begin{cases} 1 & i = j, \\ 0 & i \neq j,\end{cases} \qquad\text{for all}\qquad i,j \in \N,
	\]
	and that each element $y \in \Y$ can be represented as $y = \sum_{i \in \N} \langle y,e_i\rangle_{\Y} e_i$.}

\subsection{Detailed model assumptions}

The model \eqref{eq:inverseproblem_schematic} has to be understood in a weak sense, i.e.
\begin{equation}\label{eq:data_expandedform}
	Y_\sigma(h) = \langle Af, h\rangle_{\Y} + \sigma\xi(h), \quad h\in \Y.
\end{equation}
The error $\xi$ is a Gaussian white noise on {$\Y$}:
\begin{itemize}
	\item[(1)] If $\X$ and $\Y$ are real Hilbert spaces, we suppose that $\xi:\Y\to L^2(\Omega,P)$, for some  some probability space $(\Omega,\mathcal{A},P)$, is a linear mapping satisfying $\xi(h)\sim\mathcal{N}(0,\Vert h\Vert_\Y^2)$ and $\E\left(\xi(h)\xi(h')\right)= \langle h,h'\rangle_{\Y}$ for all $h,h'\in \Y$.
	\item[(2)] If $\X$ and $\Y$ are complex Hilbert spaces, instead we suppose that $\xi(h)\sim \mathcal{CN}(0,2\Vert h\Vert_\Y^2)$ and $\E(\xi(h)\overline{\xi(h')})= 2\langle h,h'\rangle_{\Y}$. Here $X\sim \mathcal{CN}(0,1)$ means that $X$ is distributed according to the standard complex normal distribution, i.e. $X=X_1 + i X_2$, where $X_1,X_2\iid \mathcal{N}(0,1/2)$.
\end{itemize}
We will use the notation $\langle Y_\sigma, h\rangle_\Y := Y_\sigma(h)$ for convenience.
\subsection{Notation}

For a complex number $z$, we denote its real and imaginary part by $\Re z $ and $\Im z$, respectively.\\
For two families $(a_\sigma)_{\sigma>0}$,  $(b_\sigma)_{\sigma>0}$ of non-negative real numbers we write $a_\sigma \precsim b_\sigma$ if $\lim_{\sigma\to 0} a_\sigma/b_\sigma \leq 1$, and we write $a_\sigma \succsim b_\sigma$ if $\lim_{\sigma\to 0} a_\sigma/b_\sigma \geq 1$. If $\lim_{\sigma\to 0} a_\sigma/b_\sigma = 1$, we write $a_\sigma \asymp b_\sigma$, and if $\lim_{\sigma\to 0} a_\sigma/b_\sigma = c<\infty$, we write $a_\sigma \sim b_\sigma$.

\subsection{Testing and distinguishability}\label{subsec:testing}

In the above testing problem \eqref{eq:testproblem_A}, we wish to test the \textit{hypothesis} $H_0$ against the \textit{alternative} $H_{1,\sigma}$, which means making an educated guess (based on the data) about the correctness of the hypothesis when compared to the alternative, while keeping the error of wrongly deciding against $H_0$ under control. Tests are based on \textit{test statistics}, i.e. measurable functions of the data $Y_\sigma$. We suppose that any test statistics can be expressed in terms of the Gaussian sequence $y_\sigma= (y_{\sigma,i})_{i\in \N}$ given by
\begin{equation}\label{eq:data_sequencemodel}
	y_{\sigma,i} := \langle Y_\sigma, e_i\rangle_\Y = \langle Af, e_i\rangle_\Y +\sigma\xi_{i},\quad i\in \N,
\end{equation}
and, consequently, $\xi_i\iid \mathcal{N}(0,1)$ (in the real case) or $\xi_i\iid \mathcal{CN}(0,1)$ (in the complex case) for $i\in \N$. In the following, we use the notation $Y_\sigma$ interchangeably for either the random process given by \eqref{eq:data_expandedform} or the random sequence given by \eqref{eq:data_sequencemodel}, since they are equivalent in terms of the data they provide.

A test for the testing problem \eqref{eq:testproblem_A} can now be viewed as a measurable function of the sequence $y_\sigma$ given by
\[
\phi: \mathbb{K}^{\N} \to \{0,1\}, 
\] 
where $\mathbb{K}$ is either $\R$ or $\C$. The test $\phi$ can be understood as a decision rule in the following sense: If $\phi(y_\sigma)=0$, the hypothesis is \textit{accepted}. If $\phi(y_\sigma)=1$, the hypothesis is \textit{rejected} in favor of the alternative.

If $H_0$ is true, i.e. $f=0$, but $\phi(y_\sigma)=1$, we call this a \textit{type I error} (the hypothesis is rejected although it is true). The probability to make a type I error is
\[
\alpha_\sigma(\phi):= \Prob_{0}(\phi(y_\sigma)=1),
\]
where $\Prob_0$ denotes the distribution of $y_\sigma$ given that $H_0$ is true. Likewise, the alternative might be true, but $\phi(y_\sigma)=0$. We call this a \textit{type II error} (the hypothesis is accepted although the alternative is true). Let us, for simplicity, introduce the notation $\mathcal{F}_\sigma(\mu_\sigma)=\{\delta u : u\in \mathcal{F}_\sigma, |\delta|\geq \mu_\sigma\}$. The type II error probability, given that a specific $f\in \mathcal{F}_\sigma(\mu_\sigma)$ is the true signal, is denoted as
\[
\beta_\sigma(\phi, f):= \Prob_f(\phi(y_\sigma)=0), \quad f\in\mathcal{F}_\sigma(\mu_\sigma),
\]
where $\Prob_f$ denotes the distribution of $y_\sigma$ given that $f$ is the true underlying signal. Since the alternative is -- in general -- composite, i.e. does not only consist of only one element,  the type II error probability will in general depend on the element $f$. For such composite alternatives we consider the worst case error given by the maximum type II error probability over $\mathcal{F}_\sigma(\mu_\sigma)$ for our analysis.

We say that the hypothesis $H_0$ is asymptotically \textit{distinguishable} (in the minimax sense) from the family of alternatives $(H_{1,\sigma})_{\sigma>0}$ when there exist tests for the testing problems ``$H_0$ against $H_{1,\sigma}$", $\sigma>0$, that have both small type I and small maximum type II error probabilities. We define
\[
\gamma_\sigma= \gamma_\sigma(\mu_\sigma) := \inf_{\phi \in \Phi_\sigma}\left[\alpha_\sigma(\phi)+ \sup_{f\in\mathcal{F}_\sigma(\mu_\sigma)}\beta_\sigma(\phi, f)\right],
\]
where $\Phi_\sigma$ is the set of all tests for the testing problem ``$H_0$ against $H_{1,\sigma}$". In terms of $\gamma_\sigma$ we say that $H_0$ and $H_{1,\sigma}$ are \textit{distinguishable} if $\gamma_\sigma\to 0$, as $\sigma \to 0$. If $\gamma_\sigma\to 1$, we say that they are \textit{indistinguishable}. We refer to \cite{ingster2012} for an in-depth treatment.


For prescribed families $\mathcal{F}_\sigma$, we are interested in determining the smallest possible values $\mu_\sigma$, such that $H_0$ and $H_{1,\sigma}$ are still asymptotically distinguishable, if possible. If a family $(\mu_\sigma^*)_{\sigma>0}$ exists, that satisfies
\[
\gamma_\sigma(\mu_\sigma) \to 0 \quad\text{if}\quad \mu_\sigma\succsim \mu_\sigma^*, \quad\text{and}\quad \gamma_\sigma(\mu_\sigma) \to 1 \quad\text{if}\quad \mu_\sigma\precsim \mu_\sigma^*,
\]
as $\sigma\to 0$, we call $(\mu_\sigma^*)_{\sigma>0}$the \textit{(asymptotic) minimax detection boundary}. We may say that $(\mu_\sigma^*)_{\sigma>0}$ separates \textit{detectable} and \textit{undetectable} signals. 

It is, however, not always possible to find such a sharp threshold. If the family $(\mu_\sigma^*)_{\sigma>0}$ only satisfies the weaker conditions
\[
\gamma_\sigma(\mu_\sigma) \to 0 \quad\text{if}\quad \mu_\sigma/\mu_\sigma^*\to \infty, \quad\text{and}\quad \gamma_\sigma(\mu_\sigma) \to 1 \quad\text{if}\quad \mu_\sigma/\mu_\sigma^*\to 0,
\]
we call it the \textit{separation rate} of the family of testing problems ``$H_0$ against $H_{1,\sigma}$".

\paragraph{Remark:} Although we are mostly interested in the asymptotics of the problem, we will also state non-asymptotic results, which we deem interesting.

\section{Results}\label{sec:results}

Throughout the rest of the paper, we will assume that $(u_k)_{k\in I}$ is a countable collection of functions in $\X$, and $(I_\sigma)_{\sigma>0}$ is a family of finite subsets of $I$.

\subsection{Alternatives given by finite collections of functions}\label{subsec:results_case1}

We first suppose that $\mathcal{F}_\sigma$ consists of the appropriately normalized functions $u_k$, $k\in I_\sigma$, i.e. $\mathcal{F}_{\sigma}=\left\{\Vert Au_k\Vert_\Y^{-1} u_k : k\in I_\sigma\right\}$. As above, we write $\mathcal{F}_\sigma(\mu_\sigma) = \left\{\delta \Vert Au_k\Vert_\Y^{-1} u_k : k\in I_\sigma, |\delta|\geq \mu_\sigma\right\}$, so that testing problem \eqref{eq:testproblem_A} can be written as
\begin{equation}\label{eq:testproblem_case1}
	H_0: f=0 \quad\text{against}\quad H_{1,\sigma}: f\in \mathcal{F}_\sigma(\mu_\sigma).
\end{equation}

\subsubsection{An upper bound for the detection boundary $\mu_\sigma^\ast$}

Any family of tests $(\phi_\sigma)_{\sigma>0}$ for the family of testing problems \eqref{eq:testproblem_case1} yields an upper bound for $\mu_\sigma^*$. It seems natural to choose maximum likelihood type tests as candidates, which are given by 
\begin{equation}\label{eq:MLtype_test}
	\phi_{\sigma,\alpha}(y_\sigma) = \mathbbm{1}\left\{\sup_{k\in I_\sigma} \frac{\left|\langle Y_{\sigma}, Au_k\rangle_\Y\right|}{\sigma\Vert Au_k\Vert_{\Y}}>c_{\alpha,\sigma} \right\}, \quad \sigma>0,
\end{equation}
for a given significance level $\alpha\in (0,1)$, and for appropriately chosen thresholds $c_{\alpha,\sigma}$ (which depend on whether the spaces $\X$ and $\Y$ are real or complex Hilbert spaces).

\begin{theorem}\label{thm:upperbound_case1}
	Let $N_\sigma=|I_\sigma|$ and assume that $N_\sigma\to \infty$, as $\sigma\to 0$. In addition, assume that
	\[
	\mu_\sigma \succsim (1+\varepsilon_\sigma)\sqrt{2\sigma^2\log N_\sigma},
	\]
	where $\varepsilon_\sigma\to 0$ and $\varepsilon_\sigma\sqrt{\log N_\sigma}\to \infty$ as $\sigma\to 0$. Then $\gamma_\sigma(\mu_\sigma)\to 0$ and thus, $\mu_\sigma^\ast \precsim (1+\varepsilon_\sigma)\sqrt{2\sigma^2\log N_\sigma}$.
\end{theorem}
The bound given in Theorem \ref{thm:upperbound_case1} does not depend on $A$ and it depends on set of anomalies $(u_k)_{k\in I}$ and the family of candidate indices $(I_\sigma)_{\sigma>0}$ only through the cardinality $N_\sigma$. Thus, Theorem \ref{thm:upperbound_case1} has the advantage that it is (almost) always applicable, but it might be not very well suited for specific applications. We will see examples, where the bound is essentially sharp, and an example, where it is basically useless.

\subsubsection{A lower bound for $\mu_\sigma^\ast$}

\begin{theorem}\label{thm:lowerbound_case1}
	Let
	\[
	N_\sigma^* = \sup \{\# S : S\subseteq I_\sigma,\ \Re(\langle Au_k,Au_{k'}\rangle_\Y)\leq 0\ \text{ for any two distinct } k,k'\in S \},
	\]
	and assume that $N_\sigma^*\to\infty$. In addition, assume that
	\begin{equation}\label{eq:lb}
		\mu_\sigma \precsim (1-\varepsilon_\sigma)\sqrt{2\sigma^2\log N_\sigma^*},
	\end{equation}
	where $(\varepsilon_\sigma)_{\sigma>0}$ is a family of positive real numbers such that $\varepsilon_\sigma\to 0$ and $\varepsilon_\sigma\sqrt{\log N_\sigma^*}\to \infty$ as $\sigma\to 0$. Then $\gamma_\sigma(\mu_\sigma)\to 1$ and thus, $\mu_\sigma^\ast \succsim (1-\varepsilon_\sigma)\sqrt{2\sigma^2\log N_\sigma^*}$.
\end{theorem}

This theorem can be proven by using Proposition 4.10 and Lemma 7.2 of \cite{ingster2003}. However, we will provide a self-contained and simple proof employing a weak law of large numbers for dependent random variables in section \ref{sec:proofs}.

Theorem \ref{thm:lowerbound_case1} implies that the number  $N_\sigma^*$ of negatively correlated image elements $A u_k$ is the relevant quantity which determines the difficulty of testing \eqref{eq:testproblem_A}. The actual cardinalty $N_\sigma$ plays no role in the lower bound \eqref{eq:lb}.

\subsubsection{The detection boundary}

As a consequence, we are now in position to describe the asymptotic detection boundary precisely in several situations. First, a combination of the previous theorems yields the following:
\begin{corol}\label{cor:almostsemiorth}
	Assume that $N_\sigma=|I_\sigma|\to \infty$, and let
	\[
	M_\sigma = \sup_{k\in I_\sigma}\#\{k'\in I_\sigma : \Re(\langle Au_k, Au_{k'}\rangle_\Y)> 0\},
	\]
	and assume that $M_\sigma N_\sigma^{-\varepsilon_\sigma} \to 0$ for a family $(\varepsilon_\sigma)_{\sigma>0}$ that satisfies $\varepsilon_\sigma\to 0$ and $\varepsilon_\sigma\sqrt{\log N_\sigma}\to \infty$ as $\sigma\to 0$. Then $\mu_\sigma^* \asymp \sqrt{2\sigma^2\log N_\sigma}$.
\end{corol}
In particular, Corollary \ref{cor:almostsemiorth} yields the asymptotic detection boundary, when $(Au_k)_{k \in I_\sigma}$ is orthogonal. Note that the assumptions of Corollary \ref{cor:almostsemiorth} are satisfied when $M_\sigma$ is constant as $\sigma\to 0$. This has several applications, as we will see e.g. in Section \ref{subsec:case1_examples}.

Assume now that the operator $A:\X\to\Y$ is compact and has a singular value decomposition given by orthonormal systems $(\zeta_i)_{i\in\N}$ and $(\eta_i)_{i\in\N}$ in $\X$ and $\Y$, respectively, and singular values $(s_i)_{i\in \N}$.

\begin{corol}\label{cor:case1_svd}
	Let $I=\N$ and $u_k=\zeta_k$ and $a_k=1/s_k$ for $k\in\N$, and let $(I_\sigma)_{\sigma>0}$ be any family of finite subsets of $I$, such that $N_\sigma=|I_\sigma|\to \infty$, as $\sigma\to 0$. Then $\mu_\sigma^* \asymp \sqrt{2\sigma^2\log N_\sigma}$.
\end{corol}

\paragraph{Remark:} The detection thresholds for the SVD are clearly very easy to find, and could be deduced from other known results (see \cite{ingster2012} for example). We include it here, since, as far as we know, it has not been stated explicitly before.

\subsubsection{Frame decompositions}\label{subsec:framedecomp}

We have seen that sharp detection thresholds for the SVD can easily be found, but this does (usually) not cover the situation when we are interested in local anomalies. We will thus focus on other options for anomaly systems, particularly frames, for which be briefly introduce the most important notation. Let $\mathcal{H}$ be a separable Hilbert space, and let $I$ be a countable index set. A sequence $(e_k)_{k\in I}\subseteq \mathcal{H}$ is called a frame of $\mathcal{H}$ if there exist constants $C_1,C_2>0$, such that for any $h\in\mathcal{H}$ 
\[
C_1\Vert h\Vert_\mathcal{H}^2\leq \sum_{k\in I}|\langle h,e_k\rangle_{\mathcal{H}}|^2\leq C_2\Vert h\Vert_\mathcal{H}^2.
\]

Since frames {do} not have to be orthonormal, they provide great flexibility. Theorems \ref{thm:lowerbound_case1} and \ref{thm:upperbound_case1} clearly apply to testing \eqref{eq:testproblem_A} with $u_k = e_k$, however, the fact that $(u_k)_{k \in I}$ constitutes a frame is, on its own, not enough to guarantee that we obtain a sharp detection boundary from Corollary \ref{cor:almostsemiorth}.

In the following we show how frames $(u_k)_{k \in I}$ can be constructed, for which Corollary \ref{cor:almostsemiorth} can be applied. The idea is as follows: Since the bounds for the detection threshold mostly depend on properties of the images $Au_k$ in $\Y$, we will simply start by defining a frame $(v_k)_{k\in I}$ in $\Y$ that will guarantee that the needed properties are satisfied, and then construct the corresponding frame $(u_k)_{k\in I}$ in $\X$, such that the pair $(u_k)_{k\in I}$, $(v_k)_{k\in I}$ is a decomposition of the operator $A$, and such that the assumptions of Corollary \ref{cor:almostsemiorth} are satisfied for any family of subsets $(I_\sigma)_{\sigma>0}$.

\begin{assum}\label{ass:stabilityproperty}
	\begin{itemize}
		\item[(i)] There is a dense subspace $\tilde{\Y}\subseteq \Y$ with inner product $\langle \cdot, \cdot\rangle_{\tilde{\Y}}$ and norm $\Vert \cdot\Vert_{\tilde{\Y}}$,  and constants $c_1, c_2>0$, such that 
		\begin{equation}\label{eq:stability}
			c_1\Vert x\Vert_\X \leq \Vert Ax\Vert_{\tilde{\Y}}\leq c_2\Vert x\Vert_\X,
		\end{equation}
		for all $x\in\X$.
		\item[(ii)] There is a frame $(v_k)_{k\in I}$ of $\Y$ and a sequence $(\lambda_k)_{k\in I}$ of real numbers with $\alpha_k\neq 0$, and constants $a_1,a_2>0$, such that
		\[
		a_1\Vert y\Vert_{\tilde{\Y}}^2\leq \sum_{k\in I}\lambda_k^2 |\langle y, v_k\rangle_\Y|^2\leq a_2\Vert y\Vert_{\tilde{\Y}}^2,
		\]
		for all $y\in \tilde{\Y}$.
	\end{itemize}
\end{assum}

Assumption \ref{ass:stabilityproperty} implies that $A$ as an operator from $\X$ to $\ran(A)\subseteq \tilde{\Y}$ is invertible. Now let $(v_k)_{k\in I}$ be a frame of $\ran(A)$ as in (ii). We apply the Gram-Schmidt procedure with respect to the inner product $\langle\cdot,\cdot \rangle_\Y$ to $(v_k)_{k\in I}$. This results in a sequence $(v_k^*)_{k\in I}$, which is a frame in $\ran(A)$ and which is orthogonal with respect to $\langle\cdot,\cdot \rangle_\Y$. Now we define
\[
u_k=\lambda_k A^{-1}v_k^*,
\]
for $k\in I$. The system $(u_k)_{k\in I}$ clearly yields sharp detection thresholds, as for any subset $I_\sigma \subset I$ it holds that $N_\sigma = N_\sigma^*$ by construction. Furthermore, it is a frame in $\X$, since for $x\in \X$
\[
\sum_{k\in I}|\langle x, u_k\rangle_{\X}|^2= \sum_{k\in I}\lambda_k^2|\langle (A^*)^{-1}x, v_k^*\rangle_{\Y}|^2\sim \Vert (A^*)^{-1}x\Vert_{\tilde{\Y}},
\]
and
\[
\Vert A^*\Vert_{\tilde{\Y}\to \X}^{-1}\Vert x\Vert_\X\leq \Vert (A^*)^{-1}x\Vert_{\tilde{\Y}}\leq \Vert (A^*)^{-1}\Vert_{\X\to \X}\Vert x\Vert_\X.
\]
As a consequence we obtain the following.
\begin{theorem}\label{thm:case1_framedecomp}
	Suppose that Assumption \eqref{ass:stabilityproperty} is satisfied. Then for any frame $(u_k)_{k\in I}$ of $\X$, constructed as above, and for any family of subsets of indices $(I_\sigma)_{\sigma>0}$ with $N_\sigma:=|I_\sigma|\to \infty$ as $\sigma\to 0$, we have $\mu_\sigma^*\asymp \sqrt{2\sigma^2\log N_\sigma}$.
	
\end{theorem}

\subsubsection{Examples}\label{subsec:case1_examples}

We discuss several commonly used operators and present a few typical examples of collections $(u_k)_{k\in I}$, for which the above theorems may or may not apply.

\subsubsection*{Integration}\label{subsec:integration_example}
Let $\X=\Y=L^2(\R)$ and let $A:\X\to\Y$ be the linear Fredholm integral operator given by
\[
(Af)(x) = \int_{-\infty}^x f(t)dt,\quad x\in \R, 
\]
for $f\in \X$. Suppose that $\psi$ is a (mother) wavelet in $L^2(\R)$, that satisfies $\int_\R\psi(x)dx=0$, and for which the collection $(\psi_{j,k})_{j,k\in\Z}$ given by 
\[
\psi_{j,l}(x) = 2^{j/2} \psi(2^jx-l)
\]
forms an orthogonal frame of $L^2(\R)$. For an in-depth treatment of wavelet theory, we refer to \cite{mallat} or \cite{daubechies1992ten}.

Let us suppose that the system of possible anomalies is given by this wavelet system, i.e. we consider $\{u_{(j,l)}: (j,l)\in I= \Z^2\}$ with $u_{(j,l)}= \psi_{j,l}$. Assume further that $\psi$ is compactly supported with support size $L$, which implies that for any pair of indices $(j,l)$ the number of indices $k'$, such that $\supp u_{(j,l)}\cap \supp u_{(j,l')} \neq \emptyset$ is at most $L$.

Since, in practical applications, we would not expect to be able to obtain observations on the whole plane $\R^2$, we suppose that an anomaly, if one exists, must lie within some compact subset of $\R^2$, e.g. the unit interval $[0,1]$. For some family of integers $(j_\sigma)_{\sigma>0}$ that satisfies $j_\sigma\to \infty$ as $\sigma\to 0$ we define the family $(I_\sigma)_{\sigma>0}$ of ``candidate" indices by
\begin{equation}\label{eq:integration_example_Isigma}
	I_\sigma= \left\{(j_\sigma, l) : \supp u_{(j_\sigma,l)}\subseteq [0,1]\right\}.
\end{equation}
Note that $N_\sigma\asymp 2^{j_\sigma}$. Since $\supp Af \subseteq \supp f$, it follows that for any $l$, the number of indices $l'$ such that $\supp Au_{j,l}\cap \supp Au_{j,l'}\neq \emptyset$ is bounded by $L$. Thus, the number of indices $l'$ such that  $\langle Au_{j,l}, Au_{j,l'}\rangle_\Y>0$ is also bounded by $L$. This means that $N^*_\sigma \geq N_\sigma/L$ and $M_\sigma =L$. Consequently, the conditions of Theorem \ref{cor:almostsemiorth} are satisfied, and it follows that, in this case, $\mu_\sigma^*\asymp \sqrt{2\sigma^2\log N_\sigma}$.

\subsubsection*{Periodic convolution}

Let $h:\R\to \mathbb{C}$ be a $1$-periodic and bounded function, and let $A$ be the integral operator $A:L^2([0,1])\to L^2([0,1])$ given by
\begin{equation}\label{eq:periodic_convolution}
	(Af)(x) := \int_0^1 h(u-x)f(u)du, \quad x\in [0,1].
\end{equation}
The system $(e_k)_{k\in \Z}$, where $e_k(x)=e^{-ikx}$, is a {complete orthonormal system} of $L^2([0,1])$, which consists of singular functions of $A$, since $A^*Ae_k = |\hat{h}(k)|^2e_k$. Thus, Corollary \ref{cor:case1_svd} yields the detection threshold for the detection of anomalies given by $u_k=e_k$.

Let us now try to come up with another system of possible anomalies. For the sake of simplicity, let us, from now on, only consider spaces of real-valued functions, i.e. let $\X=\Y=L^2([0,1],\R)$. Motivated by the previous example, let $\{\psi_{j,l}: j,l\in \Z\}$ be a system of compactly supported wavelets with one vanishing moment (i.e. $\int_\R x\psi_{j,l}(x)dx=0$) forming an orthonormal frame of $L^2(\R)$. We define periodic wavelets $\psi^{(per)}_{j,l}= \sum_{z\in \Z}\psi_{j,l}(\cdot+z)$ for $l=0,\ldots, 2^{j}-1$. The system $(u_{(j,l)})_{(j,l)\in I}$ given by $u_{(j,l)}=\psi^{(per)}_{j,l}$ for $I=\bigcup_{j\in\Z}\{j\}\times \{0,\ldots,2^j-1\}$ then forms an orthonormal frame of $L^2([0,1])$. 
If the function $h$ is {sufficiently} smooth, this constitutes a setting in which, for certain choices of $I_\sigma$, Theorem \ref{thm:lowerbound_case1} cannot be applied and the upper bound from Theorem \ref{thm:upperbound_case1} is basically useless, as can be seen in the following lemma.

\begin{lemma}\label{lm:deconv_wavelet_example}
	Suppose that $h:\R\to \R$ is a $1$-periodic, symmetric and continuously differentiable function, and suppose that its derivative $h'$ is Lipschitz. Let $\X=\Y=L^2([0,1],\R)$ and let $A$ be the convolution operator defined by \eqref{eq:periodic_convolution}. Let  $I =\bigcup_{j\in\Z}\{j\}\times \{0,\ldots,2^j-1\}\subseteq \R^2$ and define $u_{(j,l)}= \psi^{(per)}_{j,l}$ as above for any $(j,l)\in I$. Let $(j_\sigma)_{\sigma>0}$ be a family of integers that satisfies $j_\sigma\to \infty$ as $\sigma\to 0$ and set
	\[
	I_\sigma = \left\{(j_\sigma, l):l=0,\ldots, 2^{j_\sigma}-1\right\}. 
	\]
	for $\sigma>0$. Then $\gamma_\sigma\to 0$ if $\mu_\sigma/\sigma\to \infty$.
\end{lemma}

Intuitively, this is explained as follows. When scaled properly, the convolution of a smooth function $h$ with a wavelet with one vanishing moment on a small scale (i.e. when $j_\sigma$ is large) approximates the derivative of $h$, but shifted according to the shift parameter of the wavelet (cf. equation (6.15) of \cite{mallat}). This means that, although the support of $u_{(j_\sigma,l)}$ gets smaller when $\sigma\to 0$, the same is not true for $Au_{(j_\sigma,l)}$. In fact, it turns out, that two possible signals $Au_{(j_\sigma,l)}\Vert Au_{(j_\sigma,l)}\Vert_\Y^{-1}$ and $Au_{(j_\sigma,l')}\Vert Au_{(j_\sigma,l')}\Vert_\Y^{-1}$ will be very close (w.r.t. $\Vert\cdot\Vert_\Y$) when $l$ and $l'$ are close, and hence a test which scans over way less $k$ than in  $I_\sigma$ performs comparably well as \eqref{eq:MLtype_test}.

\subsubsection*{Radon transform}

Let us finally discuss the example of computerized tomography already mentioned in the introduction. Here, we restrict ourselves to spatial dimension $2$, in order to ease readability. We stress, however, that all subsequent results can be extended to any dimensions. Mathematically, this is modeled by the integral operator $R: L^2(\mathcal{B})\to L^2(Z, (1-t^2)^{-1/2})$, where $\mathcal{B}=\{x\in \R^2: |x|\leq 1\}$ and $Z=  [-1,1]\times[0,2\pi)$, given by
\[
(Rf)(t,\theta)= \int_\R f(t\cos \theta +s\sin \theta,t\cos \theta -s\sin \theta) ds,
\]
known as the Radon transform. The singular system of $R$ is analytically known (see \cite{natterer2001}). Let $I= \{(k,l): k\in \N_0,\, |l|\leq k,\, k+l \text{ even}\}$. We define functions $u_{(k,l)}\in L^2(\mathcal{B})$, $(k,l)\in I$ by
\[
u_{(k,l)}(x):= e^{il\varphi} r^{|l|} P_k^{(0,|l|)}(2r^2-1),\quad x=(r\cos\varphi, r\sin\varphi)\in \mathcal{B},
\]
where $P_k^{(0,|l|)}$ are the Jacobi polynomials uniquely determined by the equations $\int_0^1 t^l P_k^{(0,|l|)}P_{k'}^{(0,|l|)}=\delta_{k,k'}$. The system $(u_{(k,l)})_{(k,l)\in I}$ is a {complete orthonormal system in} $L^2(\mathcal{B})$ and, together with the appropriate {system} $(v_{(k,l)})_{(k,l)\in I}$ and constants $(\lambda_{(k,l)})_{(k,l)\in I}${,} forms the SVD of the Radon transform $R: L^2(\mathcal{B})\to L^2(Z, (1-t^2)^{-1/2})$. Thus, Corollary \ref{cor:case1_svd} yields the detection thresholds for the system $(u_{(k,l)})_{(k,l)\in I}$.

However, the discussion in Section \ref{subsec:framedecomp} gives rise to another option to choose systems of anomalies that attain the same detection boundaries. For $n\in \N$ we define the usual Sobolev space 
\[
H^\alpha(\R^n)= \left\{f\in L^2(\R^n) : \Vert f\Vert_{H^\alpha(\R^n)}<\infty\right\},
\]
where $\Vert f\Vert_{H^\alpha(\R^n)}^2:= \int_{\R^2} (1+|w|)^2)^\alpha|\hat{f}(w)|^2dw$, and set (in the notation of \cite{natterer2001})
\[
H_0^\alpha(\mathcal{B}^\circ):=\left\{f\in H^\alpha(\R^n) : \supp f\subseteq \mathcal{B}\right\}.
\]
In addition, let
\[
H^\alpha(\R\times [0,2\pi)) := \left\{f\in L^2(\R\times [0,2\pi)) : \Vert f\Vert_{H^\alpha(\R\times [0,2\pi))}<\infty\right\},
\]
where
\[
\Vert f\Vert_{H^\alpha(\R\times [0,2\pi))}^2:= \int_0^{2\pi} \Vert f(\cdot,\theta)\Vert_{H^\alpha(\R)}^2d\theta.
\]
The Radon transform is an operator from $H_0^\alpha(\mathcal{B}^\circ)$ to $H^\alpha(\R\times [0,2\pi))$ that satisfies (see Theorem 5.1 of \cite{natterer2001})
\[
C_1\Vert f\Vert_{H_0^\alpha(\mathcal{B}^\circ)}\leq \Vert Rf\Vert_{H^\alpha(\R\times [0,2\pi))}\leq C_2\Vert f\Vert_{H_0^\alpha(\mathcal{B}^\circ)}
\]
for any $f\in H_0^\alpha(\mathcal{B}^\circ)$. Thus, Theorem \ref{thm:case1_framedecomp} can be applied. The range of $R$ in $H^\alpha(\R\times [0,2\pi))$ is $\ran(R)= \left\{f\in H^\alpha(\R\times [0,2\pi)): \supp f\subseteq (-1,1)\times [0,2\pi)\right\}$. Thus, any orthonormal frame $(v_k)_{k\in I}$ of $\ran(R)$ gives rise to a frame $(u_k)_{k\in I}$ of $\X= H_0^\alpha(\mathcal{B}^\circ)$ with sharp detection boundaries given by Theorem \ref{thm:case1_framedecomp}.

\subsection{Alternatives given by the linear span of collections of anomalies}\label{subsec:results_case2}

Assume now that possibles anomalies might be linear combinations of the $u_k$, $k\in I_\sigma$. For the upcoming analysis it is necessary to assume that the $u_k$ satisfy the following.

\begin{assum}\label{assum:quasisingular}
	There is a collection $(v_k)_{k\in I}$ of functions in $\Y$, and a sequence $(\lambda_k)_{k\in I}$ of non-zero complex numbers, such that for any $f\in \X$ it holds that
	\[
	\langle Af, v_k\rangle_\Y = \lambda_k \langle f,u_k\rangle_X.
	\]
\end{assum}
Assumption \ref{assum:quasisingular} guarantees that we can present our results in terms of the $u_k$. Clearly, it is satisfied, when $u_k\in \ran(A^*)$ for all $k\in I$. In addition, if we were to assume that the collections $(u_k)_{k\in I}$ and $(v_k)_{k\in I}$ have some kind of useful structure (we may for example assume that they constitute frames of $\X$ and $\Y$, respectively, as we did in Subsection \ref{subsec:framedecomp}), then the sequence $(\lambda_k)_{k\in I}$ from Assumption \ref{assum:quasisingular} takes the role of what might be called quasi-singular values.

In this section, we suppose that $\mathcal{F}_\sigma$ consists of functions in the linear span of the functions $u_k$, $k\in I_\sigma$, namely $\mathcal{F}_\sigma=\mathcal{F}^L_{\sigma}=\left\{f\in \lspan\{u_k: k\in I_\sigma\}: \sum_{k\in I_\sigma} |\lambda_k\langle f, u_k\rangle_\X|^2 =1 \right\}$. Thus, testing problem \ref{eq:testproblem_A} becomes
\begin{equation}\label{eq:testproblem_case2}
	H_0:f=0 \quad\text{against}\quad H_{1,\sigma}: f\in \mathcal{F}_\sigma^L(\nu_\sigma),
\end{equation}
where
\[
\mathcal{F}_\sigma^L(\nu_\sigma) = \left\{f\in \lspan\{u_k: k\in I_\sigma\}: \sum_{k\in I_\sigma} |\lambda_k\langle f, u_k\rangle_\X|^2 \geq \nu_\sigma^2 \right\}.
\]
for some  family of positive real numbers $(\nu_{\sigma})_{\sigma>0}$ (we use the notation $\nu_\sigma$ instead of $\mu_\sigma$ to avoid confusion with the results from the previous section).

\subsubsection{Nonasymptotic results}
For a subset $J\subseteq I$, we define the matrix $\Xi_J$ by $(\Xi_J)_{k,k'}= \langle v_k, v_{k'}\rangle_\Y$, $k,k'\in J$, and the matrix $\tilde{\Xi}_J$ by $(\tilde{\Xi}_J)_{k,k'}= \langle\tilde{v}_k, \tilde{v}_{k'}\rangle_\Y$, $k,k'\in J$, where
\[
\tilde{v}_k := \lambda_k^{-1} Au_k,
\]
for $k\in I$. We denote the Frobenius norm of a matrix $M$ by $\Vert M\Vert_F$.

The next theorem (the non-asymptotic upper bound for the detection threshold) can not be given in terms of the minimax sum of errors $\gamma_\sigma$. Instead we define
\[
\gamma_{\sigma, \alpha}(\nu_\sigma) = \inf_{\phi \in \Phi_{\sigma,\alpha}}\left[\alpha_\sigma(\phi)+ \sup_{f\in\mathcal{F}^L_\sigma(\nu_\sigma)}\beta_\sigma(\phi, f)\right],
\]
where $\Phi_{\sigma,\alpha}$ is the set of all level $\alpha$ tests for the testing problem $``H_0 \text{ vs } H_{1,\sigma}"$. In other words, we consider the minimax sum of errors when only level $\alpha$ tests are allowed.

\begin{theorem}\label{thm:case2_upperbound_nonasymp}
	Suppose that Assumption \ref{assum:quasisingular} holds. Assume that the family of subsets $(I_\sigma)_{\sigma>0}$ is such that the matrices $\Xi_{I_\sigma}$ are positive definite for all $\sigma>0$. Then, for any $\alpha\in (0,1)$ and $\delta\in (\alpha,1)$, we have $\gamma_{\sigma,\alpha}(\nu_\sigma)\leq \delta$ if
	\[
	\nu_\sigma \geq \varepsilon d_\alpha(\delta) \sigma \sqrt{\Vert \Xi_{I_\sigma} \Vert_F},
	\]
	where $d_\alpha(\delta) = \sqrt{\log\frac{1}{\delta-\alpha}} +\left(\log\frac{1}{\alpha(\delta-\alpha)}+\sqrt{2\log\frac{1}{\delta-\alpha}}+\sqrt{2\log\frac{1}{\alpha}}\right)^{1/2}$, and $\varepsilon$ is given by $\varepsilon= 1$ if $\X$ and $\Y$ are real Hilbert spaces and $\varepsilon=\sqrt{2}$ if $\X$ and $\Y$ are complex Hilbert spaces.
\end{theorem}
It is now obvious why it is necessary to allow only tests at a prescribed level $\alpha$. Making $\alpha$ arbitrarily small would require the detection threshold to become arbitrarily large in order to keep the type II error small.

Contrary to the upper bound, the non-asymptotic lower bound for the detection threshold can be stated in terms of  $\gamma_\sigma$.

\begin{theorem}\label{thm:case2_lowerbound_nonasymp}
	Suppose that Assumption \ref{assum:quasisingular} holds, and assume that the family of subsets $(I_\sigma)_{\sigma>0}$ is such that the matrices $\tilde{\Xi}_{I_\sigma}$ are positive definite for all $\sigma>0$. Then, for any $\delta\in (0,1)$, we have $\gamma_\sigma(\nu_\sigma)\geq \delta$ if
	\[
	\nu_\sigma \leq c(\delta) \sigma \sqrt{\Vert \tilde{\Xi}_{I_\sigma}^{-1} \Vert_F},
	\]
	where $c(\delta) = \left(\log (1+(2-2\delta)^2)\right)^{1/4}$.
\end{theorem}

\paragraph{Remark 1:} The assumption that $\Xi_{I_\sigma}$ and $\tilde{\Xi}_{I_\sigma}$, respectively, are positive definite (and consequently invertible, since they are Hermitian) is a technical necessity. However, it is also intuitively justified, because it prevents certain ``unreasonable" choices of $I_\sigma$ (for example any subset $I_\sigma$ such that $(u_k)_{k\in I_\sigma}$ is linearly dependent).

\paragraph{Remark 2:} Note that it can be easily seen that, if we consider the set 
\[
\left\{f\in \lspan\{u_k: k\in I_\sigma\}: \sum_{k\in I_\sigma} |\langle f, u_k\rangle_\X|^2 \geq \nu_\sigma^2 \right\},
\]
instead of $\mathcal{F}_\sigma^{L}(\nu_\sigma)$, then we would obtain the same bounds as above with $\Xi_{I_\sigma}$ replaced by the matrix $\Lambda_{I_\sigma}$, which is given by $(\Lambda_{I_\sigma})_{k,k'} = (\lambda_k \overline{\lambda_{k'}})^{-1} \langle v_k,v_{k'}\rangle_\Y$, and $\tilde{\Xi}_{I_\sigma}$ replaced by the matrix $\tilde{\Lambda}_{I_\sigma}$ given by $(\tilde{\Lambda}_{I_\sigma})_{k,k'} = \langle Au_k,Au_{k'}\rangle_\Y$ It follows, that our results are compatible with the results obtained in \cite{laurentloubes2012}, where the above testing problem was considered when the system $(u_k, v_k, \lambda_k)_{k\in I}$ is given by the SVD of $A$.

\subsubsection{Asymptotic results}

The asymptotic results for this section can now be easily deduced from the previous theorems.

\begin{corol}\label{cor:case2_asymp} Suppose that the assumptions of Theorems \ref{thm:case2_upperbound_nonasymp} and \ref{thm:case2_lowerbound_nonasymp} hold.
	\begin{itemize}
		\item[(1)]  $H_0$ and $H_{1,\sigma}$ are asymptotically distinguishable if
		\[
		\frac{\nu_{\sigma}}{\sigma \sqrt{\Vert \Xi_{I_\sigma}\Vert_F}}\to \infty.
		\]
		\item[(2)] $H_0$ and $H_{1,\sigma}$ are asymptotically indistinguishable if
		\[
		\frac{\nu_{\sigma}}{\sigma \sqrt{\Vert \tilde{\Xi}^{-1}_{I_\sigma}\Vert_F}}\to 0.
		\]
	\end{itemize}
\end{corol}

Until now, we have allowed the $u_k$ to just be any functions we might be interested in detecting. However, we are able to refine our results when we assume that $(v_k)_{k\in I}$ and $(\tilde{v}_k)_{k\in I}$ are ``well-behaved". We call a sequence of functions $(h_i)_{i\in\N}$ from some Hilbert space $\mathcal{H}$ a \textit{Riesz sequence}, if there exist constants $C_1, C_2>0$ such that
\[
C_1\sum_{i\in \N}|a_i|^2 \leq \Vert \sum_{i\in \N}a_i h_i\Vert_\mathcal{H}^2\leq C_2\sum_{i\in \N}|a_i|^2,
\]
for any sequence $(a_i)_{i\in\N}\in \ell^2$. Two sequences $(h_i)_{i\in\N}$ and $(h'_i)_{i\in\N}$ are called \textit{biorthogonal} if
\[
\langle h_i, h'_j\rangle_\mathcal{H}=\delta_{i,j},
\]
where $\delta_{i,j}$ is the Kronecker symbol.

\begin{assum}\label{assum:case2_frames}
	The collections $(v_k)_{k\in I}$ and $(\tilde{v}_k)_{k\in I}$ are biorthogonal Riesz sequences.
\end{assum}
We acknowledge that Assumption \ref{assum:case2_frames} is restrictive. We will discuss non-trivial situations in which it is satisfied below. We collect the implications of Assumption \ref{assum:case2_frames} in the following lemma.

\begin{lemma}\label{lm:rate_for_constant_norm}
	Suppose that Assumptions \ref{assum:quasisingular} and \ref{assum:case2_frames} hold, and let $(I_\sigma)_{\sigma>0}$ be an arbitrary family of subsets of $I$. Then the following statements hold.
	\begin{itemize}
		\item[(1)] For any $\sigma>0$, the matrices $\Xi_{I_\sigma}$ and $\tilde{\Xi}_{I_\sigma}$ are positive definite.
		\item[(2)] There are constants $c_1,c_2>0$ such that $c_1\Vert \Xi_{I_\sigma}\Vert_F\leq \Vert \tilde{\Xi}^{-1}_{I_\sigma}\Vert_F\leq c_2\Vert \Xi_{I_\sigma}\Vert_F$.
		\item[(3)] $\Vert \Xi_{I_\sigma}\Vert_F\sim N_\sigma^{1/2}$ (and consequently, also $\Vert \tilde{\Xi}^{-1}_{I_\sigma}\Vert_F\sim N_\sigma^{1/2}$), as $\sigma\to 0$, where $N_\sigma= |I_\sigma|$.
	\end{itemize}
\end{lemma}

Thus, if all conditions of Lemma \ref{lm:rate_for_constant_norm} are satisfied, it follows from Corollary \ref{cor:case2_asymp} that the \textit{separation rate} of the family of testing problems ``$H_0$ against $H_{1,\sigma}$" is given by $\nu_\sigma^*\sim \sigma N_\sigma^{1/4}$.

\paragraph{Remark:} Note that this result is not surprising. The separation rate corresponds to the rate (in terms of the euclidean norm) of detecting an $n$-dimensional signal $\theta\in\R^n\setminus\{0\}$ from observations given by $X=\theta + \sigma Z$, where $Z\sim \mathcal{N}(0,\id_n)$. Furthermore, it has been shown previously (cf. \cite{laurentloubes2012}) that the same holds, when $(u_k,v_k,\lambda_k)_{k\in I}$ constitute the SVD of the operator $A$. Thus, the above results yield a generalization of the known theory.

\subsubsection{Examples}\label{subsec:case2_examples}
It is clear that, when the system $(u_k, v_k, \lambda_k)_{k\in I}$ is given by the SVD of the operator $A:\X\to \Y$, then all of the above theory can be applied. Since this was the subject of \cite{laurentloubes2012}, we will omit a discussion of this example here.

\subsubsection*{Examples based on the wavelet-vaguelette decomposition}
Suppose that $(u_k)_{k\in I}$ is a system of orthogonal wavelets in $\X$. If chosen appropriately (for a complete discussion, see \cite{donoho1995wvd}), it follows that for certain operators $A:\X\to \Y$ , there exist non-zero numbers $(\lambda_k)_{k\in I}$, such that the systems $(v_{k})_{k\in I}$ and $(\tilde{v}_{k})_{k\in I}$ of functions in $\Y$ given by
\[
A^*v_{k} = \lambda_k u_{k}, \quad \tilde{v}_{k} = \lambda_k^{-1} A u_{k}
\]
form biorthogonal Riesz sequences in $\Y$ (see Theorem 2 of \cite{donoho1995wvd})). Clearly, Assumption \ref{assum:quasisingular} is satisfied in this case.

We immediately see that this would yield nice examples of the theory developed in this section. We will discuss a few situations, in which such a construction is possible, below.

\subsubsection*{Integration}\label{subsec:case2_int_example}

Consider the setting of example \ref{subsec:integration_example}, i.e. $u_{(j,l)}=\psi_{j,k}$ for $(j,l)\in I=\Z^2$ and for some wavelet $\psi$. Suppose that the wavelet $\psi$ is continuously differentiable. In this case, the WVD is particularly simple. Let
\[
v_{(j,l)}(x)= -2^{j/2}\psi'(2^jx-l), \quad \tilde{v}_{(j,l)}(x)= 2^{j/2}\psi^{(-1)}(2^jx-l),
\]
with $\lambda_{(j,l)}=2^{-j}$. Then it follows from \cite{donoho1995wvd} that the systems $(v_{(j,l)})_{(j,l)\in I}$ and $(\tilde{v}_{(j,l)})_{(j,l)\in I}$ form biorthogonal Riesz sequences in $L^2(\R)$. Thus, we can apply Lemma \ref{lm:rate_for_constant_norm} to obtain $\Vert \Xi_{I_\sigma}\Vert_F\sim N_\sigma^{1/2}$, and thus, $\nu_\sigma^*\sim \sigma N_\sigma^{1/4}$ for any family $(I_\sigma)_{\sigma>0}$ of ``candidate" indices.

\subsubsection*{Periodic Convolution}
Let $\X= \Y= L^2(\mathcal{S}^1)$, where $\mathcal{S}^1$ is the unit circle. In other words we consider square-integrable $1$-periodic functions on $[0,1]$. Let the operator $A$ be given by $(Af)(x) = \int_0^1 h(x-u)f(u)du$, $x\in [0,1]$ for some $1$-periodic function $h$. Let $(u_{j,k})_{(j,k)\in \Z^2}$ be a basis of periodic Meyer wavelets, each with Fourier coefficients $u_{j,k,m}$, $m\in\Z$, and let $(I_\sigma)_{\sigma>0}$ be any family of finite subsets of $I=\Z^2$. Let $h_m$, $m\in \Z$ be the Fourier coefficients of $h$. It was shown in Appendix B of \cite{johnstone2004} that, if $h_m= C|m|^{-a}$, for some $a>0$, the collections $(v_{j,k})_{(j,k)\in\Z^2}$ and  $(\tilde{v}_{j,k})_{(j,k)\in\Z^2}$ given by
\[
v_{j,k}(x) = \sum_{m\in \Z} \frac{\lambda_{j,k} u_{j,k,m}}{h_m}e^{imx},\quad \tilde{v}_{j,k}(x) = \sum_{m\in \Z} \frac{h_m u_{j,k,m}}{\lambda_{j,k}}  e^{imx},
\]
where the quasi-singular values are given by $\lambda_{j,k}= 2^{-ja}C$, yield biorthogonal Riesz sequences.

\subsubsection*{Radon transform}
We start by introducing two-dimensional wavelet systems. Let $\psi$ be a (one-dimensional) wavelet and $\varphi$ its corresponding scaling function. We assume that they are of compact support and at least two times continuously differentiable. Define the two-dimensional functions
\[
\eta^{0}(x) := \varphi(x_1)\theta(x_2),\ \  \eta^{1}(x) := \psi(x_1)\varphi(x_2),\ \ \eta^{2}(x) := \varphi(x_1)\psi(x_2),\ \ \eta^{3}(x) := \psi(x_1)\psi(x_2),
\]
where $x=(x_1,x_2)\in\R^2$. The function $\eta^{0}$ is the two-dimensional scaling function, and the functions $\eta^{\varepsilon}$, $\varepsilon\in\{1,2,3\}$, are the two-dimensional wavelets. They are scaled and translated as usual, i.e.
\[
\eta^{\varepsilon}_{j,l}(x) := 2^j \eta^{\varepsilon}(2^j x -l), \quad j\in \Z,\ l\in \Z^2,\ \varepsilon\in \{0,1,2,3\}.
\]
The collection of all translated and scaled wavelets (without the scaling functions, i.e. excluding $\varepsilon=0$) is a {complete orthonormal system} of $L^2(\R^2)$ (see for example Theorem 7.24 of \cite{mallat}). Using the projection theorem (see Theorem 1.1 of \cite{natterer2001}), it was shown in \cite{donoho1995wvd} that for $\varepsilon\in\{1,2,3\}$ and for $f\in L^2(\R^2)\cap \dom(R)$, we have
\[
\int_0^{\pi}\int_{-\infty}^{\infty} (Rf)(t,\theta) \overline{(R\omega^{\varepsilon}_{j,k})(t,\theta)} dtd\theta = \int_{\R^2} f(x)\overline{\eta^{\varepsilon}_{j,k}(x)} dx 
\]
where $\omega^{\varepsilon}_{j,k}$ is defined through its Fourier transform $\widerhat{\omega^{\varepsilon}_{j,k}}(x) = \frac{1}{2\pi} |x| \widerhat{\eta^{\varepsilon}_{j,k}}(x)$. In other words,
\begin{equation}\label{eq:radon_quasisingula_rrelation}
	\langle Rf, R\omega^{\varepsilon}_{j,k}\rangle_{L^2(\R\times [0,2\pi))}  = \langle f, \eta^{\varepsilon}_{j,k}\rangle_{L^2(\R^2)}.
\end{equation}
An in-homogeneous wavelet basis of $L^2(\R)$ can be constructed as follows: For some $j_0\in \Z$, we consider the collection of functions
\[
\left\{\eta^0_{j_0,k} : k\in \Z^2\right\} \cup \left\{\eta^\varepsilon_{j,k} : \varepsilon\in \{1,2,3\}, k\in\Z^2, j\geq j_0\right\}.
\]
Note that for any $j_0\in \Z$, $k_0\in \Z^2$ we can write
\[
\eta^{0}_{j_0,k_0} =\sum_{\varepsilon\in\{1,2,3\}} \sum_{k\in \Z^2} \sum_{ -\infty<j\leq j_0} c^{\varepsilon}_{j_0,k_0,j,k} \eta^{\varepsilon}_{j,k},
\]
and thus, it follows from the linearity of all the above operations that the relation \eqref{eq:radon_quasisingula_rrelation} is also true for $\varepsilon=0$ with $w^{0}_{j,k}$ defined accordingly.

With practical applications (where it is an unreasonable assumption that observations on all of the plane $\R^2$ can be made) in mind, we assume that signals, if existent, lie within a compact set, e.g. the unit ball $\mathcal{B}= \{x\in \R^2 : \Vert x\Vert_2\leq 1\}$. We consider the Radon transform as an operator $R:\X\to \Y$, where $\X=L^2(B)$, and $\Y= L^2(Z)$ with $Z=[-1,1]\times [0,\pi)$. Note that, contrary to the example in Section \ref{subsec:case1_examples}, the space $\Y$ is equipped with the norm $\Vert \cdot\Vert_{\Y}$ given by $\Vert f\Vert_\Y = \int_{-1}^1 \int_0^{2\pi} |f(t,\theta)|^2d\theta dt$. (Note that the operator $\R:\X\to\Y$ is well-defined and bounded since $\Vert f\Vert_{L^2(Z)}\leq \Vert f\Vert_{L^2(Z, (1-t)^{-1/2})}$ for any $f\in L^2(Z)$.)

We devise a ``wavelet-type" frame of $L^2(\mathcal{B})$ as follows. We choose $j_0$ large enough, such that the area of $\supp(\eta^\varepsilon_{j_0,k})$ is small compared to the area of the unit ball $\mathcal{B}$. Now let
\[
I = \left\{(j_0,k,0) : \supp(\eta^0_{j_0,k})\cap \mathcal{B} \neq \varnothing\right\} \cup \left\{(j,k,\varepsilon) : j\geq j_0,\ \varepsilon\in \{1,2,3\},\ \supp(\eta^\varepsilon_{j,k})\cap \mathcal{B} \neq \varnothing\right\}, 
\]
and finally, define $u_{(j,k,\varepsilon)} = \eta^\varepsilon_{j,k}|_{\mathcal{B}}$ for $(j,k,\varepsilon)\in I$. The collection $(u_{(j,k,\varepsilon)})_{(j,k,\varepsilon)\in I}$ forms a frame of $L^2(\mathcal{B})$, since, for any $f\in L^2(\R^2)$ supported in $\mathcal{B}$, we have
\[
\langle f, u_{(j,k,\varepsilon)}\rangle_{L^2(\mathcal{B})} = \langle f, \eta^\varepsilon_{j,k}\rangle_{L^2(\R^2)}.
\]
Note that $\Vert R\omega^{\varepsilon}_{(j,k)}\Vert_{L^2(\R\times [0,2\pi))}\sim 2^{j/2}$ (again, see \cite{donoho1995wvd}). Thus, if we let $v_{(j,k,\varepsilon)} = 2^{-j/2}R\omega^{\varepsilon}_{(j,k)}|_{Z}$, we obtain
\[
\langle Rf, v_{(j,k,\varepsilon)}\rangle_{L^2(Z)} = 2^{-j/2}\langle Rf, R\omega^{\varepsilon}_{(j,k)} \rangle_{L^2(\R\times[0,2\pi))}= 2^{-j/2}\langle f, \eta^{\varepsilon}_{(j,k)}\rangle_{L^2(\R^2)} =  2^{-j/2}\langle f, u_{(j,k,\varepsilon)}\rangle_{L^2(\mathcal{B})},
\]
for any $f\in L^2(\mathcal{B})$ (which we extended to $L^2(\R^2)$ by setting $f(x)=0$ whenever $x\not\in \mathcal{B}$). Thus, Assumption \ref{assum:quasisingular} is satisfied for the set $(u_{(j,k,\varepsilon)})_{(j,k,\varepsilon)\in I}$ with $v_{(j,k,\varepsilon)}$ defined as above and $\lambda_{(j,k,\varepsilon)}= 2^{-j/2}$.

It follows that Theorems \ref{thm:case2_upperbound_nonasymp} and \ref{thm:case2_lowerbound_nonasymp} are applicable for the collection $(u_{(j,k,\varepsilon)})_{(j,k,\varepsilon)\in I}$ and yield non-asymptotic results for appropriate choices of $I_\sigma$. However, note that the $u_{(j,k,\varepsilon)}$ are not necessarily orthonormal. 

It follows from Lemma 4 (and the discussion leading up to it) of \cite{donoho1995wvd} that, if  $\psi$ has at least 4 vanishing moments and is at least 4 times continuously differentiable, then the collections $(2^{-j/2}R\omega^{\varepsilon}_{j,k})_{j\in\Z, k\in \Z^2,\varepsilon\in\{1,2,3\}}$ and $(2^{j/2}R\eta^{\varepsilon}_{j,k})_{j\in\Z, k\in \Z^2,\varepsilon\in\{1,2,3\}}$ are Riesz sequences.

If we suppose that all subsets $I_\sigma$ are chosen such that all $\eta^\varepsilon_{j,k}$ lie completely within $\mathcal{B}$, i.e. $\supp(\eta^\varepsilon_{j,k})\subseteq \mathcal{B}$ for any $(j,k,\varepsilon)\in I_\sigma$ for all $\sigma>0$, then $\tilde{v}_{(j,k,\varepsilon)}= 2^{j/2}R\eta^{\varepsilon}_{j,k}$, and it follows from the same arguments as in the proof of Lemma \ref{lm:rate_for_constant_norm} that
\[
\Vert \tilde{\Xi}_{I_\sigma}^{-1}\Vert_F\geq C\Vert \Omega_{I_\sigma}\Vert_F\geq C'\sqrt{N_\sigma},
\]
where $\Omega_{I_\sigma}$ is given by $(\Omega_{I_\sigma})_{k,k'} = \langle 2^{-j/2}R\omega^{\varepsilon}_{j,k}, 2^{-j/2}R\omega^{\varepsilon}_{j,k'}\rangle_{L^2(\R\times [0,2\pi))}$. On the other hand, since $(v_k)_{k\in I}$ is a frame of $L^2([-1,1]\times[0,\pi))$, it follows (as in the proof of Lemma \ref{lm:rate_for_constant_norm}) that 
\[
\Vert \Xi_{I_\sigma}\Vert_F^2\leq C'' \sum_{k\in I_\sigma} \Vert v_k\Vert_{L^2([-1,1]\times[0,\pi)} \leq C'' \sum_{k\in I_\sigma} \Vert 2^{-j/2}R\omega^{\varepsilon}_{j,k} \Vert_{L^2(\R\times[0,\pi)} \leq C''' N_\sigma.
\]
Thus, $\nu_\sigma^*\sim \sigma N_\sigma^{1/4}$.

\section{Simulation study}\label{sec:simulations}
\subsection*{A note on discretization}

\renewcommand{\vol}{\mathrm{vol}}

For convenience, we will from now on assume that $\mathbb{K}=\R$. Assume that $\Y=L^2(D,\R)$ for some $D\subseteq \R^d$, $d\in\N$. For a finite subset $S\subseteq D$, we define the evaluation function $e_S$ by $e_S: h\mapsto (h(s))_{s\in S}$. Now let $A_S = e_S\circ A$. Clearly, $A_S: \X\to \R^{n}$, where $n=|S|$, is a bounded linear operator. We equip $\R^{n}$ with the inner product $\langle \cdot,\cdot\rangle_n$ (and the corresponding norm $\Vert \cdot\Vert_n$) given by $\langle x, x'\rangle_n^2 = \frac{\vol(D)}{n} \sum_{s\in S}x_s x'_s$, for $x=(x_s)_{s\in S}, x'=(x'_s)_{s\in S}\in \R^{n}$, and thereby make $\R^{n}$ a Hilbert space. Here, $\vol(D)$ denotes the ($d$-dimensional) volume of $D$. Now suppose that we observe data $Y_{\sigma,S}$ on $\R^{n}$ given by
\begin{equation}\label{eq:data_discretized}
	Y_{\sigma,S} = A_S(f)+ \sigma \sqrt{\frac{n}{\vol(D)}} \xi,
\end{equation}
where $\xi=(\xi_s)_{s\in S}\sim \mathcal{N}(0,\id_n)$. Since, for any $x,x'\in \R^{n}$, we have
\[
\left\langle \sqrt{\frac{n}{\vol(D)}}\xi, x\right\rangle_n \sim \mathcal{N}(0, \Vert x\Vert_n)\quad \text{and}\quad \E\left(\left\langle \sqrt{\frac{n}{\vol(D)}} \xi, x\right\rangle_n \left\langle \sqrt{\frac{n}{\vol(D)}} \xi, x'\right\rangle_n\right) = \langle x,x'\rangle_n,
\]
it follows that our results are valid for this discretized model. Note that asymptotic results still refer to $\sigma$ becoming small (and not $n$ becoming large). If $S$ is chosen appropriately, $\langle e_S(h), e_S(h')\rangle_n$ can be viewed as an approximation of $\langle h,h'\rangle_\Y$ for $h,h'\in \Y$, and, consequently, the upper and lower bounds derived from \eqref{eq:data_discretized} can be viewed as an approximation of the upper and lower bounds derived from the data \eqref{eq:inverseproblem_schematic}.

Finally, note that testing $``f=0"$ against $``f\in \mathcal{F}_\sigma(\mu_\sigma)"$ based on the data $Y_{\sigma,n}$ is equivalent to testing $``f=0"$ against $``f\in \mathcal{F}_\sigma\left(\frac{\mu_\sigma\vol(D)}{n}\right)"$ based on $X_{\sigma,S}$ given by
\[
X_{\sigma,S}:= A_S(f)+ \sigma \sqrt{\frac{\vol(D)}{n}} \xi.
\]

\subsection*{Integration}
We consider the example from section \ref{subsec:integration_example}, discretized as above with $S=\{\frac{i}{n}: i=0,\ldots,n-1\}$ for $n=2^{15}$. The wavelet system $(\psi_{j,k})_{j,k\in\Z}$ consists of Daubechies (db6) wavelets. See Figure \ref{fig:case1_int} for the results of the simulation study. Note that the displayed results are only approximations in two senses: First, we used the test $\phi_{\sigma,\alpha}= \phi_{\sigma,\alpha}^{ML}$ from \eqref{eq:MLtype_test}, which may not necessarily be the optimal test, and second, we approximate $\sup_{f\in\mathcal{F}_\sigma(\delta)}\beta_\sigma(\phi_{\sigma,\alpha}^{ML}, f)$ by  $\beta^*_\sigma(\delta) := N_\sigma^{-1}\sum_{k\in I_\sigma}\beta_\sigma(\phi_{\sigma,\alpha}^{ML}, \delta \Vert Au_k\Vert_\Y^{-1} u_k)$, i.e. the mean type II error over all possible anomalies of minimal ``amplitude", since it is, in general, not clear, which $k\in I_\sigma$ will maximize the type II error. 

Note that the proof of Theorem \ref{thm:upperbound_case1} yields a non-asymptotic upper bound for the detection threshold: We have $\sup_{f\in\mathcal{F}_\sigma(\delta)}\beta_\sigma(\phi_{\sigma,\alpha}^{ML}, f)\leq \alpha$ when $\delta\geq \sigma(c_{\alpha,\sigma}+z_{1-\alpha})$, where $z_{1-\alpha}$ is the $(1-\alpha)$-quantile of the standard Gaussian distribution.

Next, we consider the example from section \ref{subsec:case2_int_example}. Everything is as above, except for a few differences: We consider alternatives given by linear combinations of $(u_k)_{k\in I_\sigma}$ as in Section \ref{subsec:results_case2}, we use the test $\phi_{\sigma,\alpha}= \phi_{\sigma,\alpha}^{\chi^2}$ given by \eqref{eq:chisquaredtest}, and $\beta^*_\sigma(\delta)$ is given by $\E_\pi \beta_\sigma(\phi_{\sigma,\alpha}^{ML}, \delta f)$, where $\pi$ is the uniform distribution on $\mathcal{F}_\sigma^L$. The results of this study are displayed in Figure \ref{fig:case2_int}.

\subsection*{Radon transform}
The setting for our simulation study for the Radon transform is inspired by the discussion in Section \ref{subsec:case2_examples}. We consider the Radon transform as an operator
\[
R:L^2(\mathcal{B}_{\frac{1}{\sqrt{2}}}, \R) \to L^2\left(\left[-\frac{1}{\sqrt{2}},\frac{1}{\sqrt{2}}\right]\times [0,\pi), \R\right),
\]
where $\mathcal{B}_{\frac{1}{\sqrt{2}}}= \left\{x\in \R^2 : \Vert x\Vert_2\leq \frac{1}{\sqrt{2}}\right\}$ is the ball that contains the unit square $[-1/2,1/2]^2$. Let $\{\eta_{j,l}^\varepsilon: j\in \Z, l\in \Z^2, \varepsilon\in \{1,2,3\}\}$ be the two-dimensional wavelet system (consisting of Daubechies (db4) wavelets) from Section \ref{subsec:case2_examples}, define $u_{(j,l,\varepsilon)}=\eta_{j,l}^\varepsilon$ for $(j,l,\varepsilon)\in I$ with
\[
I= \left\{(j,l,\varepsilon) : \supp \eta_{j,l}^\varepsilon\subseteq[-1/2,1/2]^2\right\},\]
and let $I_\sigma= \left\{(j_\sigma,l,\varepsilon) : \supp \eta_{j_\sigma,l}^\varepsilon\subseteq[-1/2,1/2]^2\right\}$, for family $(j_\sigma)_{\sigma>0}$ of natural numbers. We consider discretized data of the form \eqref{eq:data_discretized} with $S=\left\{\left(-\frac{1}{2}+\frac{i_1}{1024},\frac{i_2}{360}\pi\right): i_1=0,\ldots,1023, i_2=0,\ldots,359\right\}$. As above, we use the test $\phi_{\sigma,\alpha}= \phi_{\sigma,\alpha}^{\chi^2}$ given by \eqref{eq:chisquaredtest}, and $\beta^*_\sigma(\delta)=\E_\pi \beta_\sigma(\phi_{\sigma,\alpha}^{ML}, \delta f)$, where $\pi$ is the uniform distribution on $\mathcal{F}_\sigma^L$. The results of this study are displayed in Figure \ref{fig:case2_radon}. 

The example in Figure \ref{fig:radonintro} also comes from this setting: The parameters were $j_\sigma=5$, $\sigma=15$, and $\delta= 264$. By Theorem \ref{thm:case2_upperbound_nonasymp}, the distorted image can be distinguished from the reference image with type I and type II error both at most $0.05$ by the test $\phi_{8,0.05}^{\chi^2}$.

\begin{figure}[H]
	\begin{subfigure}{.45\textwidth}
		\centering
		\includegraphics[width=\linewidth]{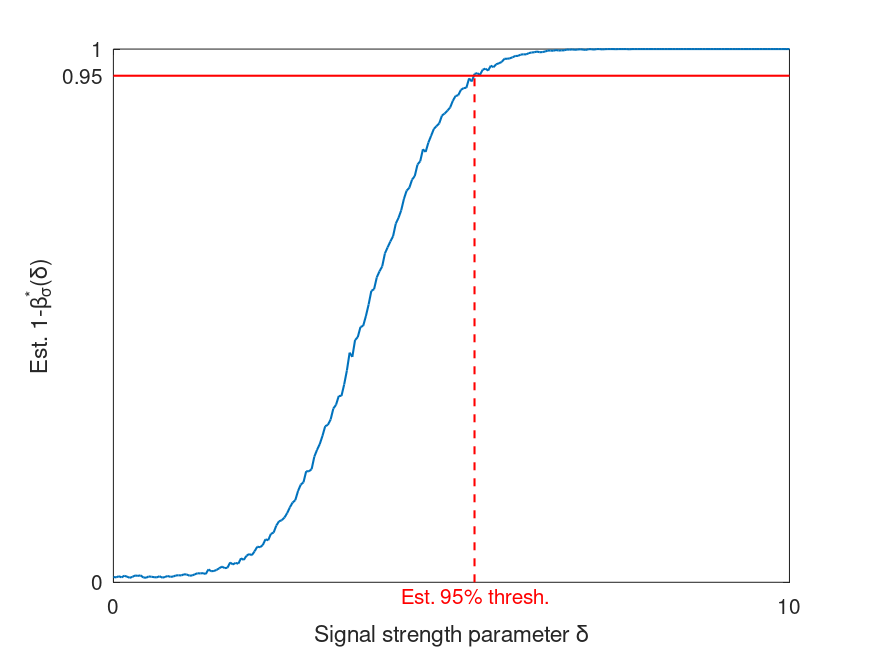}
	\end{subfigure}
	\begin{subfigure}{.45\textwidth}
		\centering
		\includegraphics[width=\linewidth]{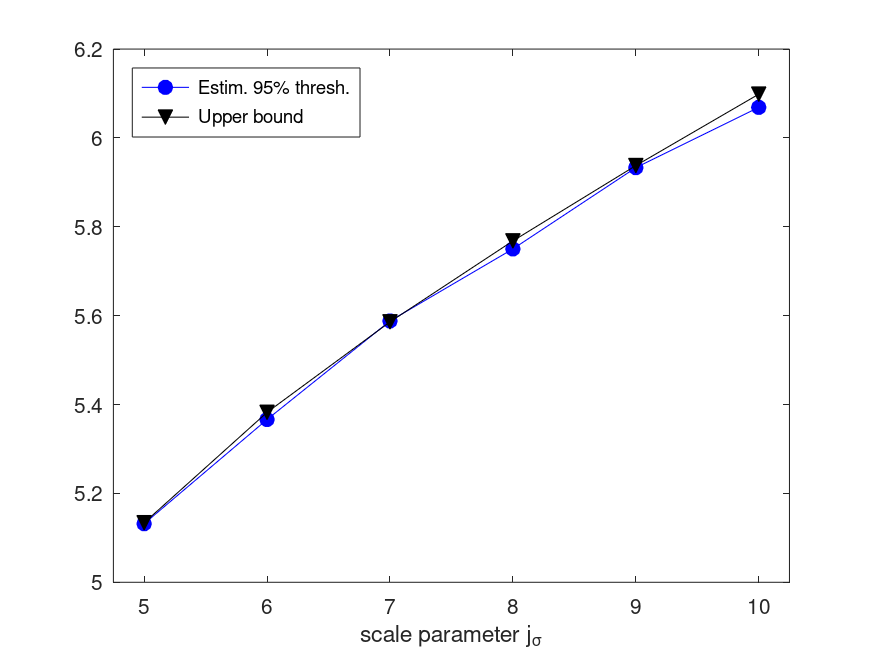}
	\end{subfigure}
	\caption{\textbf{Left:} Estimation of $1-\beta^*_\sigma(\delta)$ for $\delta$ between $0$ and $10$, for $j_\sigma=6$, $\sigma=1$ and  $\alpha=0.05$. For each value of $\delta$, $5000$ tests have been performed. The results suggest that the power achieves $95\%$ for $\delta\approx 5.3414$. \textbf{Right:} Estimated values of $\delta$ for which the power achieves $95\%$ for $j_\sigma\in \{5,\ldots,10\}$ and $\sigma=1$ compared with the upper bound derived from the proof of Theorem \ref{thm:upperbound_case1}.}\label{fig:case1_int}
\end{figure}

\begin{figure}[H]
	\begin{subfigure}{.45\textwidth}
		\centering
		\includegraphics[width=\linewidth]{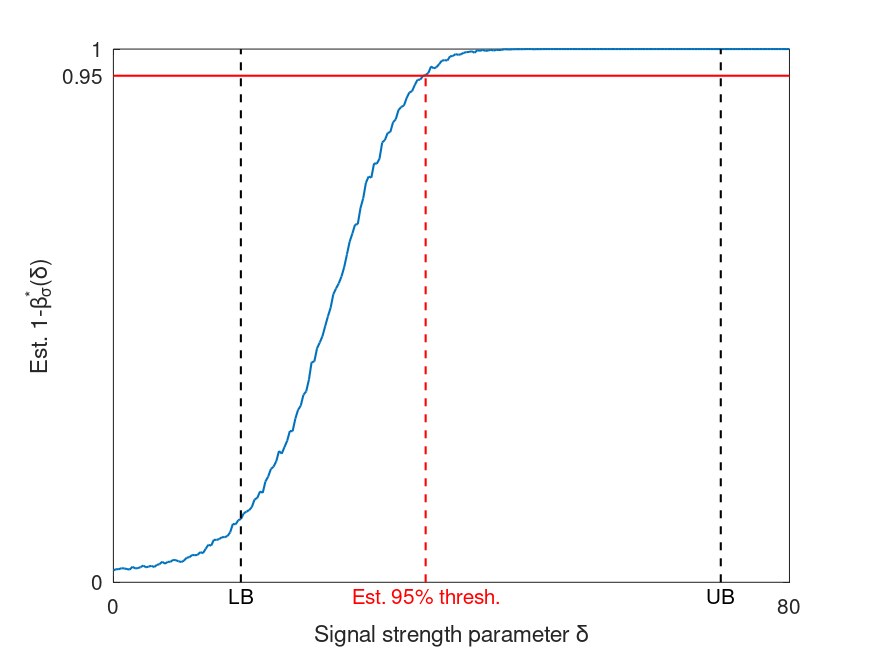}
	\end{subfigure}
	\begin{subfigure}{.45\textwidth}
		\centering
		\includegraphics[width=\linewidth]{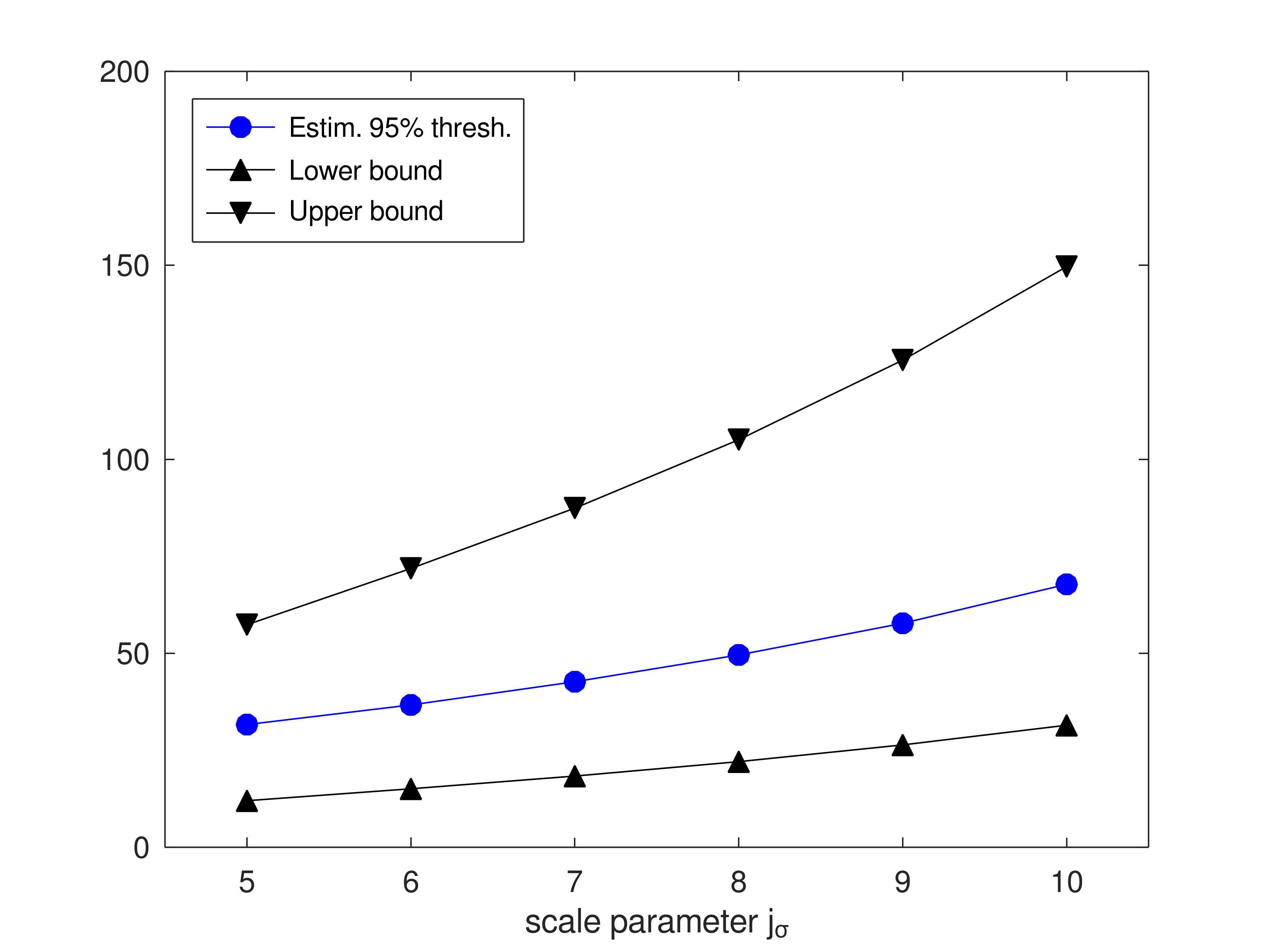}
	\end{subfigure}
	\caption{\textbf{Left:} Estimation of $1-\beta^*_\sigma(\delta)$ for $\delta$ between $0$ and $80$, for $j_\sigma=6$, $\sigma=1$ and  $\alpha=0.05$, compared with the upper bound (UB) from \ref{thm:case2_upperbound_nonasymp} and the lower bound (LB) from Theorem \ref{thm:case2_lowerbound_nonasymp}. For each value of $\delta$, $5000$ tests have been performed. The results suggest that the power achieves $95\%$ for  $\delta\approx 36.948$. \textbf{Right:} Estimated values of $\delta$ for which the power achieves $95\%$ for $j_\sigma\in \{5,\ldots,10\}$ and $\sigma=1$ compared with the upper bound from \ref{thm:case2_upperbound_nonasymp} and the lower bound from Theorem \ref{thm:case2_lowerbound_nonasymp}.}\label{fig:case2_int}
\end{figure}

\begin{figure}[H]
	\begin{subfigure}{.45\textwidth}
		\centering
		\includegraphics[width=\linewidth]{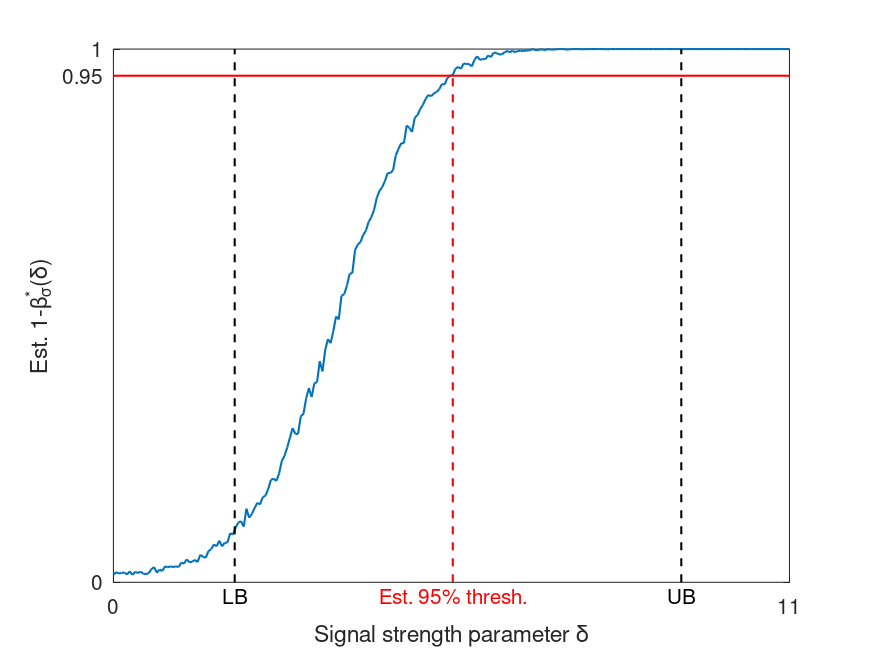}
	\end{subfigure}
	\begin{subfigure}{.45\textwidth}
		\centering
		\includegraphics[width=\linewidth]{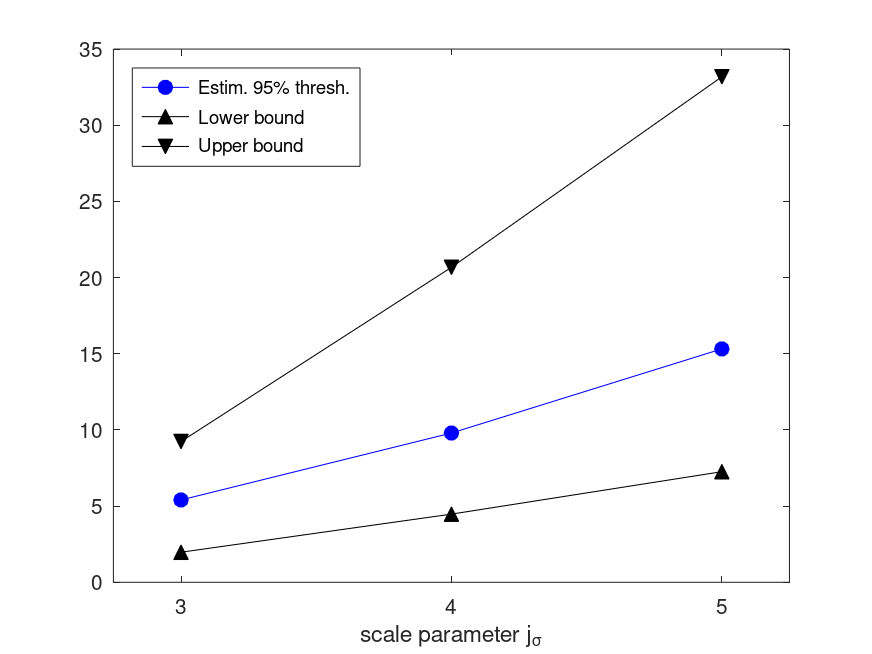}
	\end{subfigure}
	\caption{\textbf{Left:} Estimation of $1-\beta^*_\sigma(\delta)$ for $\delta$ between $0$ and $11$, for $j_\sigma=3$, $\sigma=1$ and  $\alpha=0.05$, compared with the upper bound (UB) from \ref{thm:case2_upperbound_nonasymp} and the lower bound (LB) from Theorem \ref{thm:case2_lowerbound_nonasymp}. For each value of $\delta$, $5000$ tests have been performed. The power achieves $95\%$ for $\delta\approx 5.5221$. \textbf{Right:} Estimated values of $\delta$ for which the power achieves $95\%$ for $j_\sigma\in \{5,\ldots,10\}$ and $\sigma=1$ compared with the upper bound from \ref{thm:case2_upperbound_nonasymp} and the lower bound from Theorem \ref{thm:case2_lowerbound_nonasymp}.}\label{fig:case2_radon}
\end{figure}

These simulations seem to affirm our theoretical results. Note that the thresholds displayed in Figure \ref{fig:case1_int} are very large compared to the asymptotic detection boundary from Corollary \ref{cor:almostsemiorth}. This is due to the logarithmic growth of the detection boundary.

\section{Discussion}
In this paper, we have considered statistical hypothesis testing in inverse problems with localized alternatives. This can be used to determine whether an unknown object that deviates from a reference object, can be distinguished from that reference or not.

More precisely, we first considered alternatives given by finitely many elements (e.g. chosen from a dictionary), and under additional restrictions on the structure of this system we were able to derive the (asymptotic) detection boundary. Those results are illustrated along examples such as integration, convolution, and the Radon transform. Afterwards, we have moved to more complex alternatives allowing for linear combinations of elements from the dictionary. In this case, we were still able to derive the minimax separation rate even under weaker assumptions on the structure of the system. This has been illustrated again for the above-mentioned and in simulations.

The results in this study offer several point of contact for further research. For practical purposes, the design of more (computationally and statistically) efficient multiple tests is on demand, which is beyond the scope of this paper. It would be interesting to see which methods can be used to efficiently test a reference object against hypotheses which consist e.g. of wavelets on different scales, which is a setting, for which the assumptions of Corollary \ref{cor:almostsemiorth} are not satisfied in general. Another interesting question is the detection boundary in case of sparse alternatives (similar to \cite{laurentloubes2012}), which we have not discussed here. 

\section*{Acknowledgements}
M. P. acknowledges support from the RTG 2088. A. M. acknowledges support of the DFG Cluster of Excellence 2067 ``Multiscale Bioimaging: From Molecular Machines to Networks of Excitable Cells". F. W. is supported by the DFG via grant WE 6204/4-1. The authors would like to thank Markus Haltmeier and Miguel del \'Alamo for insightful comments. Furthermore we would like to express our gratitude to two anonymous referees, whose comments helped to improve the presentation of the paper.

\section{Proofs}\label{sec:proofs}

\subsection{Proof for section \ref{subsec:results_case1}}

\subsubsection{Proof of the upper bound}

\begin{proof}[Proof of Theorem \ref{thm:upperbound_case1}] We treat the two cases (whether $\X$ and $\Y$ are real or complex spaces) separately.
	
	\paragraph{$\X$ and $\Y$ are real Hilbert spaces.}
	Any test for the testing problem \eqref{eq:testproblem_A} yields an upper bound for $\gamma_\sigma$, and, thus, also an upper bound for $\mu_\sigma^*$. Our upper bound is based on a particularly simple family of likelihood ratio type tests given by \eqref{eq:MLtype_test} with thresholds given by
	\[
	c_{\alpha,\sigma}= \sqrt{2\log \frac{N_\sigma}{\alpha}}.
	\]
	We show that for any $\sigma>0$ and any $\alpha\in (0,1)$, the test $\phi_{\alpha,\sigma}$ has level $\alpha$ and its asymptotic type II error vanishes for the testing problem \eqref{eq:testproblem_A} if $\mu_\sigma \succsim (1+\varepsilon_\sigma)\sqrt{2\sigma^2\log N_\sigma}$. This would then imply that $\gamma_\sigma\precsim \alpha + o(1)$, which will immediately prove the theorem, since $\alpha$ was arbitrary.
	
	Setting $f_k=\frac{u_{k'}}{\Vert Au_{k'}\Vert_\Y}$, we have
	\[
	\frac{\langle Y_\sigma, Au_k\rangle_\Y}{\sigma \Vert Au_k\Vert_\Y} \overset{H_0}{\sim} \mathcal{N}(0,1) \quad\text{and}\quad  \frac{\langle Y_\sigma, Au_k\rangle_\Y}{\sigma \Vert Au_k\Vert_\Y} \overset{f=\delta f_k}{\sim} \mathcal{N}\left(\frac{\delta \langle Au_{k'}, Au_k\rangle_{\Y}}{\sigma\Vert Au_{k'}\Vert_\Y\Vert Au_k\Vert_\Y},1\right)
	\]
	Using the union bound and a concentration inequality for the normal distribution we find
	\begin{align*}
		\Prob_{H_0}\left(\phi_{\alpha,\sigma}(y_\sigma)= 1\right) &= \Prob_{H_0}\left(\sup_{k\in I_\sigma}  \frac{\left|\langle Y_\sigma, Au_k\rangle_\Y\right|}{\sigma \Vert Au_k\Vert_\Y}>c_{\alpha,\sigma}\right)\\
		&\leq N_\sigma \Prob\left( \left|Z\right|> c_{\alpha,\sigma}\right) =N_\sigma \exp\left(-\frac{1}{2}c_{\alpha,\sigma}^2\right)=\alpha,
	\end{align*}
	for some $Z\sim\mathcal{N}(0,1)$. Thus, $\phi_{\alpha,\sigma}$ is indeed a level $\alpha$ test. Next, we show that the maximal type II error of $\phi_{\alpha,\sigma}$ vanishes. For some {identically distributed (but not necessarily independent)} random variables $Z_{\sigma,k}{\sim}\mathcal{N}(0,1)$, $k\in I_\sigma$, we find {by the union bound, that}
	\begin{align*}
		\sup_{k\in I_\sigma}\sup_{|\delta| \geq \mu_\sigma} \Prob_{\delta f_k}\left(\phi_{\alpha,\sigma}(y_\sigma)=0\right)&= \sup_{k\in I_\sigma}\sup_{|\delta| \geq \mu_\sigma} \Prob_{\delta f_k}\left(\sup_{k'\in I_\sigma} \frac{\left|\langle Y_\sigma, Au_{k'}\rangle_\Y\right|}{\sigma \Vert Au_{k'}\Vert_\Y}\leq c_{\alpha,\sigma}\right)\\
		&\leq \sup_{k\in I_\sigma}\sup_{|\delta| \geq \mu_\sigma} \Prob\left(\sup_{k'\in I_\sigma} \left(Z_{\sigma,k'} + \frac{\delta \langle Au_{k'}, Au_k\rangle_{\Y}}{\sigma \Vert Au_{k'}\Vert_\Y\Vert Au_k\Vert_\Y}\right) \leq c_{\alpha,\sigma}\right)\\
		&\leq \sup_{k\in I_\sigma}\sup_{|\delta| \geq \mu_\sigma} \Prob\left( Z_{\sigma,k} + \frac{\delta}{\sigma} \leq c_{\alpha,\sigma}\right)\\
		&=  \Prob\left(Z + \frac{\mu_\sigma}{\sigma} \leq c_{\alpha,\sigma}\right)\\
		&= \Prob\left(Z\leq \sqrt{2\log N_\sigma} + \sqrt{2\log(1/\alpha)}- (1+\varepsilon_\sigma)\sqrt{2\log N_\sigma}\right) \to 0,
	\end{align*}
	since $\varepsilon_\sigma\sqrt{\log N_\sigma}\to\infty$.
	
	\paragraph{$\X$ and $\Y$ are complex Hilbert spaces.}
	The idea of the proof is the same as above. We again use is the test given by \eqref{eq:MLtype_test} with thresholds
	\[
	c_{\alpha,\sigma}= \sqrt{1+2\sqrt{\log\frac{N_\sigma}{\alpha}}+2\log\frac{N_\sigma}{\alpha}}.
	\]
	Setting $f_{k'}=\frac{u_{k'}}{\Vert Au_{k'}\Vert_\Y}$ as above, we have
	\[
	\frac{\langle Y_\sigma, Au_k\rangle_\Y}{\sigma \Vert Au_k\Vert_\Y} \overset{H_0}{\sim} \mathcal{CN}(0,2) \quad\text{and}\quad  \frac{\langle Y_\sigma, Au_k\rangle_\Y}{\sigma \Vert Au_k\Vert_\Y} \overset{f=\delta f_{k'}}{\sim} \mathcal{CN}\left(\frac{\delta \langle Au_{k'}, Au_k\rangle_{\Y}}{\sigma\Vert Au_{k'}\Vert_\Y\Vert Au_k\Vert_\Y},2\right)
	\]
	We first show that $\Phi_\alpha$ is a level $\alpha$ test. For some $Z\sim\mathcal{CN}(0,1)$ we find, using the union bound, that
	\[
	\Prob_{H_0}\left(\phi_{\alpha,\sigma}(y_\sigma)= 1\right) = \Prob_{H_0}\left(\sup_{k\in I_\sigma}  \frac{\left|\langle Y_\sigma, Au_k\rangle_\Y\right|}{\sigma \Vert Au_k\Vert_\Y}>c_{\alpha,\sigma}\right)\leq N_\sigma \Prob\left( \left|Z\right|^2> c_{\alpha,\sigma}^2\right).
	\]
	Note that $|Z|^2= \Re(Z)^2+\Im(Z)^2\sim \chi_2^2$. It follows from Lemma 1 of \cite{laurentmassart2000} that
	\[
	N_\sigma \Prob\left(|Z|^2> 1+2\sqrt{\log\frac{N_\sigma}{\alpha}}+2\log\frac{N_\sigma}{\alpha}\right) \leq N_\sigma\exp\left(-\log\frac{N_\sigma}{\alpha}\right)= \alpha.
	\]
	Thus, $\phi_{\alpha,\sigma}$ is indeed a level $\alpha$ test. As above, we must now show that the maximal type II error of $\phi_{\alpha,\sigma}$ vanishes. For some {identically distributed (but not necessarily independent)} random variables $Z_{\sigma,k}{\sim}\mathcal{CN}(0,2)$, $k\in I_\sigma$ and $Z_\sigma\sim\mathcal{CN}(0,2)$, we find {by the union bound, that}
	\begin{align*}
		\sup_{k\in I_\sigma}\sup_{|\delta| \geq \mu_\sigma} \Prob_{\delta f_k}\left(\phi_{\alpha,\sigma}(y_\sigma)=0\right)&= \sup_{k\in I_\sigma}\sup_{|\delta| \geq \mu_\sigma} \Prob_{\delta f_k}\left(\sup_{k'\in I_\sigma} \frac{\left|\langle Y_\sigma, Au_{k'}\rangle_\Y\right|}{\sigma \Vert Au_{k'}\Vert_\Y}\leq c_{\alpha,\sigma}\right)\\
		&\leq \sup_{k\in I_\sigma}\sup_{|\delta| \geq \mu_\sigma} \Prob\left(\sup_{k'\in I_\sigma} \left| Z_{\sigma,k'} + \frac{\delta \langle Au_{k'}, Au_k\rangle_{\Y}}{\sigma \Vert Au_{k'}\Vert_\Y\Vert Au_k\Vert_\Y}\right|^2 \leq c_{\alpha,\sigma}^2\right)\\
		&\leq \sup_{k\in I_\sigma}\sup_{|\delta| \geq \mu_\sigma} \Prob\left(\left| Z_{\sigma,k} + \frac{\delta}{\sigma}\right|^2 \leq c_{\alpha,\sigma}^2\right) =  \sup_{|\delta| \geq \mu_\sigma} \Prob\left(\left| Z_{\sigma} + \frac{\delta}{\sigma}\right|^2 \leq c_{\alpha,\sigma}^2\right)
	\end{align*}
	We have
	\[
	\left| Z_{\sigma} + \frac{\delta}{\sigma}\right|^2 = |Z_{\sigma}|^2 + 2\Re\left(\frac{\bar{\delta}}{\sigma} Z_{\sigma}\right)+\frac{|\delta|^2}{\sigma^2}\geq \frac{2|\delta|}{\sigma} Z_{\sigma}' +\frac{|\delta|^2}{\sigma^2},
	\]
	for some $Z_{\sigma}'\sim \mathcal{N}(0,1)$. It follows that
	\begin{align*}
		\sup_{|\delta| \geq \mu_\sigma} \Prob\left(\left| Z_{\sigma} + \frac{\delta}{\sigma}\right|^2 \leq c_{\alpha,\sigma}^2\right) &\leq \Prob\left(\frac{2\mu_\sigma}{\sigma} Z_{\sigma}' +\frac{\mu_\sigma^2}{\sigma^2} \leq 1+2\sqrt{\log\frac{N_\sigma}{\alpha}}+2\log\frac{N_\sigma}{\alpha}\right)\\
		&\leq \Prob\left(Z_{\sigma}'\leq \sqrt{\frac{1}{2}\log N_\sigma} -(1+\varepsilon_\sigma)\sqrt{\frac{1}{2}\log N_\sigma} +O(1) \right) \to 0,
	\end{align*}
	since $\varepsilon_\sigma\sqrt{\log N_\sigma}\to\infty$.
\end{proof}

\subsubsection{Proof of the lower bound}

We suppose that $\X$ and $\Y$ are complex Hilbert spaces. The proof for the real case is analogous. In fact, the proof of Theorem \ref{thm:lowerbound_case1} can in principle be derived from Proposition 4.10 and Lemma 7.2 from \cite{ingster2003} with just a few adjustments.

\begin{proof}[Proof of Theorem \ref{thm:lowerbound_case1}]

	Let $I_\sigma^*$ be the largest subset of $I_\sigma$ such that $\Re(\langle Au_k, Au_{k'}\rangle_\Y)\leq 0$ for any distinct $k,k'\in I_\sigma^*$. Recall that $|I_\sigma^*|=N_\sigma^*$. Recall that
	\[
	y_{\sigma,i} := \langle Y_\sigma, e_i\rangle_\Y = \langle Af, e_i\rangle_\Y + \sigma\xi_i,
	\]
	where $\xi_i\iid \mathcal{CN}(0,2)$ for $i\in \N$. This means that under $H_0$ the random sequence
	\[
	\tilde{y}_\sigma := \left(\frac{1}{\sigma}\Re(y_{\sigma,1}),\frac{1}{\sigma}\Im(y_{\sigma,1}), \frac{1}{\sigma}\Re(y_{\sigma,2}), \frac{1}{\sigma}\Im(y_{\sigma,2}), \ldots\right)
	\]
	is a sequence of i.i.d. $\mathcal{N}(0,1)$-distributed random variables. Any test statistic may be expressed in terms of the Gaussian sequence $\tilde{y}_\sigma$. Hence, any test may be expressed as a function of $\tilde{y}_\sigma$.
	
	\paragraph{Bayesian alternative.} For $k\in I_\sigma$ we define $f_k:=\frac{\mu_\sigma u_k}{\Vert Au_k\Vert_\Y}$, and let $\pi_\sigma$ be the prior distribution on the alternative set $\mathcal{F}_\sigma(\mu_\sigma)$ given by $\pi_\sigma = \frac{1}{N_\sigma^*}\sum_{k\in I_\sigma^*} \delta_{f_k}$. The idea is to bound the maximal type II error probability from below by the mean (in terms of $\pi_\sigma$) type II error probability as follows:
	\[
	\gamma_\sigma = \inf_{\phi \in \Phi_\sigma}\left[\alpha_\sigma(\phi)+ \sup_{f\in\mathcal{F}_\sigma(\mu_\sigma)}\beta_\sigma(\phi, f)\right] \geq \inf_{\phi \in \Phi_\sigma}\left[\alpha_\sigma(\phi)+ \frac{1}{N_\sigma^*}\sum_{k\in I_\sigma^*} \beta_\sigma(\phi, f_k)\right].
	\]
	We may say that it suffices to analyze the ``simpler" testing problem
	\[
	H_0: f\equiv 0 \quad\text{vs}\quad H_{\pi,\sigma}: f\sim \pi_\sigma,
	\]
	in terms of its mean type II error, instead of \eqref{eq:testproblem_A} in the minimax sense. We have (cf. \cite{ingster2003}, chapter 2)
	\[
	\inf_{\phi \in \Phi_\sigma}\left[\alpha_\sigma(\phi)+ \frac{1}{N_\sigma^*}\sum_{k\in I_\sigma^*} \beta_\sigma(\phi, f_k)\right] =1-\frac{1}{2}\E_0\left|\E_{\pi_\sigma} \frac{d\Prob_{f}}{d\Prob_0}(\tilde{y}_\sigma)-1\right|,
	\]
	where $\E_{\pi_\sigma} \frac{d\Prob_{f}}{d\Prob_0}(\tilde{y}_\sigma)= \frac{1}{N_\sigma^*}\sum_{k\in I_\sigma^*}\frac{d\Prob_{f_k}}{d\Prob_0}(\tilde{y}_\sigma)$. We see that, in order to show indistinguishability, it suffices to show that 
	\begin{equation}\label{eq:aux_case1_lb_lln}
		\E_0\left|\E_{\pi_\sigma} \frac{d\Prob_{f}}{d\Prob_0}(\tilde{y}_\sigma)-1\right| \to 0.
	\end{equation}
	We denote the distribution of $\tilde{y}_{\sigma}$ on $\R^\N$ under $H_0$ by $\Prob_0$ (this is of course the standard Gaussian distribution on $\R^\N$). The Cameron-Martin space corresponding to $(\R^\N, \Prob_0)$ is $\mathcal{H}=\ell^2$ with norm $\Vert \cdot\Vert_\mathcal{H}^2=\Vert \cdot\Vert_2^2$ (see for example Example 4.1 of \cite{lifshits2012}). If $f=f_k$, then $\tilde{y}_\sigma$ has distribution $\Prob_{f_k}$ defined by $\Prob_{f_k}(\cdot)= \Prob_0(\cdot - h_k)$, where $h_k\in \R^\N$ is given by
	\[ 
	h_{k,2j-1}= \frac{\mu_\sigma}{\sigma \Vert Au_k\Vert_\Y}\Re(\langle Au_k, e_j\rangle_\Y),\quad
	h_{k,2j}= \frac{\mu_\sigma}{\sigma\Vert Au_k\Vert_\Y}\Im(\langle Au_k, e_j\rangle_\Y), \quad j\in \N.
	\]
	It follows that
	\[
	\Vert h_k\Vert_\mathcal{H}^2= \sum_{i\in\N}\frac{\mu_\sigma^2}{\sigma^2\Vert Au_k\Vert_\Y^2}|\langle Au_k,e_i\rangle_\Y|^2 = \frac{\mu_\sigma^2}{\sigma^2},
	\]
	which also shows that $h\in \mathcal{H}$.
	Thus, by the Cameron-Martin theorem (see Theorem 5.1 of \cite{lifshits2012}), 
	\begin{align*}
		\frac{d\Prob_{f_k}}{d\Prob_0}(\tilde{y}_\sigma) &= \exp\left(\sum_{i\in \N} h_{k,i} \tilde{y}_{\sigma,i} -\frac{\mu_\sigma^2}{2\sigma^2} \right)\\
		&= \exp\left(\frac{\mu_\sigma}{\sigma\Vert Au_k\Vert_\Y}\sum_{j\in \N} (\Re(\langle Au_k, e_j\rangle_\Y) \tilde{y}_{\sigma,2j-1} +\Im(\langle Au_k, e_j\rangle_\Y) \tilde{y}_{\sigma,2j}) -\frac{\mu_\sigma^2}{2\sigma^2} \right)\overset{H_0}{=} \exp\left(X_{\sigma,k} \right),
	\end{align*}
	where $X_{\sigma,k}\sim\mathcal{N}\left(-\frac{\mu_\sigma^2}{2\sigma^2},\frac{\mu_\sigma^2}{\sigma^2}\right)$. Note that the distribution of $X_{\sigma,k}$ does not depend on $k$. However, the collection $\{X_{\sigma,k} : k\in I_\sigma\}$ is, in general, not independent. We have
	\[
	\E_{\pi_\sigma} \frac{d\Prob_{f}}{d\Prob_0}(Y_\sigma)=\frac{1}{N_{\sigma}^*} \sum_{k\in I^*_{\sigma}}\exp(X_{\sigma,k}).
	\]
	In order to show that \eqref{eq:aux_case1_lb_lln} holds, we employ a weak law of large numbers, namely Theorem 3.2 from \cite{wang2014}. However, note that similar ideas have been used in \cite{ingster2012}.
	
	\paragraph{A weak law of large numbers.}  Let $(\sigma_m)_{m\in \N}$ be a sequence of positive real numbers, such that $\sigma_m\searrow 0$ as $m\to\infty$.
	Consider the triangular array of random variables $\{\exp(X_{\sigma_m,k}) : m\in\N, k\in I_{\sigma_m}^* \}$. Note that 
	\[
	\Cov\left(\exp(X_{\sigma_m,k}),\exp(X_{\sigma_m,k'})\right)= \exp\left(\frac{\mu_{\sigma_m}^2\Re(\langle Au_k, Au_{k'}\rangle_\Y)}{\sigma_m^2\Vert Au_k\Vert_\Y\Vert Au_{k'}\Vert_\Y}\right)-1\leq 0,
	\]
	for any $m$ and any two distinct $k,k'\in I^*_{\sigma_m}$. Let $(\kappa_m)_{m\in \N}$ be another sequence of real numbers given by $\kappa_m=(N_{\sigma_m}^*)^{(1+\varepsilon_{\sigma_m})(1-\varepsilon_{\sigma_m})^2}$. Then $\kappa_m\to\infty$ as $m\to \infty$ and
	\[
	(N_{\sigma_m}^*)^{-1}\kappa_m = (N_{\sigma_m}^*)^{-\varepsilon_{\sigma_m} +O(\varepsilon_{\sigma_m}^2)} \to 0.
	\]
	It follows that
	\begin{align*}
		\frac{1}{N_{\sigma_m}^*} \sum_{k\in I^*_{\sigma_m}} \E_0\left[\left(\exp(X_{\sigma_m,k})\right)\mathbbm{1}(\exp(X_{\sigma_m,k})>\kappa_m)\right] &= \Prob\left[\mathcal{N}(0,1)\leq \frac{\mu_{\sigma_m}}{2\sigma_m} - \frac{{\sigma_m} \log \kappa_m}{\mu_{\sigma_m}}\right]\\
		&\leq  \Prob\left[\mathcal{N}(0,1)\leq  -\varepsilon_{\sigma_m}(1+\varepsilon_{\sigma_m})\sqrt{\frac{1}{2} \log N_{\sigma_m}^*}\right],
	\end{align*}
	which vanishes as $m\to\infty$, since $\varepsilon_\sigma\sqrt{2 \log N_\sigma^*}\to\infty$. We can now employ Theorem 3.2 from \cite{wang2014}, which immediately yields that
	\[
	\E_0\left|\E_{\pi_\sigma} \frac{d\Prob_{f}}{d\Prob_0}(Y_\sigma)-1\right|= \E_0\left|\frac{1}{N_{\sigma_m}^*} \sum_{k\in I^*_{\sigma_m}}\exp(X_{\sigma_m,k})-1\right|\to 0,
	\]
	as $m\to\infty$.
	
\end{proof}

\subsubsection{Remaining proofs}
\begin{proof}[Proof of Corollary \ref{cor:almostsemiorth}]
	For $k\in I_\sigma$, let
	\[
	S_\sigma(k) = \{k'\in I_\sigma : \Re(\langle Au_k, Au_{k'}\rangle_\Y)\geq 0\}.
	\]
	We construct a subset $I_\sigma'$ of $I_\sigma$ iteratively as follows. We choose $k_1\in I_\sigma$ arbitrarily, then choose $k_2\in I_\sigma\setminus S_\sigma(k_1)$ arbitrarily, then choose $k_3\in I_\sigma\setminus (S_\sigma(k_1)\cup S_\sigma(k_2))$ arbitrarily, and continue until $S_\sigma(k_1)\cup S_\sigma(k_2)\cup\ldots=I_\sigma$. Then set $I_\sigma'=\{k_1,k_2,k_3,\ldots\}$. Since, by assumption, $|S_\sigma(k)|\leq M_\sigma$ for any $k\in I_\sigma$, it follows that $|I_\sigma'|\geq N_\sigma/M_\sigma\succsim N_\sigma^{1+\varepsilon_\sigma}$. Since the set $I_\sigma^*$ can be constructed as above, Theorem \ref{thm:lowerbound_case1} yields
	\[
	\mu_\sigma^*\succsim (1-\varepsilon_\sigma)\sqrt{1+\varepsilon_\sigma}\sqrt{2\log N_\sigma}.
	\]
	Thus,
	\[
	\sqrt{1-\varepsilon_\sigma-\varepsilon_\sigma^2+\varepsilon_\sigma^3}\precsim \frac{\mu_\sigma^*}{\sqrt{2\log N_\sigma}}\precsim 1+\varepsilon_\sigma,
	\]
	and the claim follows.
\end{proof}

\begin{proof}[Proof of Corollary \ref{cor:case1_svd}]
	Since $A\zeta_k = s_k\eta_k$ for all $k$, and the system $(\eta_k)_{k\in \N}$ is orthonormal, this follows immediately from Corollary \ref{cor:almostsemiorth}.
\end{proof}

\begin{proof}[Proof of Lemma \ref{lm:deconv_wavelet_example}]
	Since, by assumption, $\mu_\sigma/\sigma\to \infty$, we can choose a family of positive integers $(n_\sigma)_{\sigma>0}$, such that $n_\sigma\to\infty$ as $\sigma\to 0$ and
	\[
	\frac{\mu_\sigma}{\sigma}-\sqrt{2\log n_\sigma} \to \infty.
	\]
	For $m\in\{0,\ldots, n_\sigma-1\}$ let
	\[
	w_{\sigma,m} = 2^{j_\sigma/2}\psi\left(2^{j_\sigma}\left(\cdot -m/n_\sigma\right)\right),
	\]
	and let $w_{\sigma,m}^{(per)}=\sum_{z\in\Z}w_{\sigma,m}(\cdot+z)$.
	For $\alpha\in (0,1)$ consider the test $\phi_{\alpha, \sigma}(y_\sigma)=\mathbbm{1}\{T_\sigma>c_{\alpha,\sigma}\}$ with threshold $c_{\alpha, \sigma}= \sqrt{2\log \left(n_\sigma/\alpha\right)}$ and test statistic
	\[
	T_\sigma= \sup_{0\leq m\leq n_\sigma-1}  \frac{\left|\langle Y_\sigma, Aw_{\sigma,m}^{(per)}\rangle_\Y \right|}{\sigma \Vert Aw_{\sigma,m}^{(per)}\Vert_\Y}.
	\]
	It is easy to see that $\phi_{\alpha, \sigma}$ is a level $\alpha$ test. Let $f_l=\frac{\psi_{j_\sigma,l}^{(per)}}{\Vert A\psi_{j_\sigma,l}^{(per)}\Vert_\Y}$. For $l\in \{0,\ldots, 2^{j_\sigma}-1\}$ we define $m^*(l) = \argmin\{|2^{-j_\sigma}l-n_\sigma^{-1} m|: m\in \{0,\ldots, n_\sigma-1\}\}$. As in the previous proofs we find
	\[
	\sup_{0\leq l\leq 2^{j_\sigma}-1} \sup_{\delta \geq \mu_\sigma} \Prob_{\delta f_l}\left(\phi_{\alpha,\sigma}(y_\sigma)=0\right)\leq \Prob\left(Z \leq c_{\alpha,\sigma}- \frac{\mu_\sigma}{\sigma} \inf_{0\leq l\leq 2^{j_\sigma}-1}\frac{\left\langle A\psi_{j_\sigma, l}^{(per)}, Aw_{j_\sigma, m^*(l)}^{(per)}\right\rangle_\Y}{\Vert A\psi_{j_\sigma,l}^{(per)}\Vert_\Y\Vert Aw_{\sigma,m^*(l)}^{(per)}\Vert_\Y} \right),
	\]
	for some $Z\sim \mathcal{N}(0,1)$. It remains to show that
	\[
	\inf_{0\leq l\leq 2^{j_\sigma}-1}\frac{\left\langle A\psi_{j_\sigma, l}^{(per)}, Aw_{j_\sigma, m^*(l)}^{(per)}\right\rangle_\Y}{\Vert A\psi_{j_\sigma,l}^{(per)}\Vert_\Y\Vert Aw_{\sigma,m^*(l)}^{(per)}\Vert_\Y}\to 1,
	\]
	as $\sigma\to 0$. Note that, due to periodicity, for any $l\in\{0,\ldots, 2^{j_\sigma}-1\}$ and $m\in \{0,\ldots, n_\sigma-1\}$, 
	\begin{align*}
		\Vert A\psi_{j_\sigma,l}^{(per)}\Vert_\Y^2=\Vert Aw_{\sigma,m}^{(per)}\Vert_\Y^2=\Vert A\psi_{j_\sigma,0}^{(per)}\Vert_\Y^2 &= \int_0^1 \left(\int_0^1 h(u-x)\psi_{j_\sigma,0}^{(per)}(u)du\right)^2dx\\
		&= \int_0^1 \left(\int_{-\infty}^\infty h(u-x)2^{j_\sigma/2}\psi(2^{j_\sigma}u)du\right)^2dx.
	\end{align*}
	It follows from equation (6.15) of \cite{mallat}, that for any $x\in\R$,
	\[
	\lim_{\sigma\to 0} 2^{3j_\sigma/2}\int_{-\infty}^\infty h(u-x) 2^{j_\sigma/2}\psi\left(2^{j_\sigma}u\right)du = Ch'(x),
	\]
	for some constant $C>0$. Due to Theorem 6.2 of \cite{mallat}, there exist an integrable function $\theta$, such that $-\frac{d}{dx}\theta(x) = \psi(x)$. Since $\psi$ has bounded support, $\theta$ has bounded support as well. Thus, we find by substituting and integrating by parts that
	\begin{align*}
		2^{3j_\sigma/2}\int_{-\infty}^\infty h(u-x) 2^{j_\sigma/2}\psi\left(2^{j_\sigma}u\right)du &= 2^{3j_\sigma/2}\int_{-\infty}^\infty h(2^{-j_\sigma}v-x) 2^{-j_\sigma/2}\psi\left(v\right)dv\\
		&= -\int_{-\infty}^\infty h'(2^{-j_\sigma}v-x) \theta\left(v\right)dv.
	\end{align*}
	Since, by assumption, $h'$ is Lipschitz and periodic, it follows that it is bounded. Thus, for any $x\in\R$,
	\[
	\left|2^{3j_\sigma/2}\int_{-\infty}^\infty h(u-x) 2^{j_\sigma/2}\psi\left(2^{j_\sigma}u\right)du\right|\leq C'\int_{-\infty}^\infty |\theta\left(v\right)|dv.
	\]
	It follows from the dominated convergence theorem that
	\[
	2^{3j_\sigma/2}\Vert A\psi_{j_\sigma,0}^{(per)}\Vert_\Y \to C\Vert h'\Vert_\Y,
	\]
	as $\sigma\to \infty$. We have
	\[
	\frac{\left\langle A\psi_{j_\sigma, l}^{(per)}, Aw_{j_\sigma, m^*(l)}^{(per)}\right\rangle_\Y}{\Vert A\psi_{j_\sigma,l}^{(per)}\Vert_\Y\Vert Aw_{\sigma,m^*(l)}^{(per)}\Vert_\Y} = 1- \frac{\Vert A\psi_{j_\sigma, l}^{(per)} - Aw_{j_\sigma, m^*(l)}^{(per)}\Vert_\Y^2}{\Vert A\psi_{j_\sigma,0}^{(per)}\Vert_\Y^2}.
	\]
	We substitute twice, integrate by parts, use that $h'$ is Lipschitz and that $\left|2^{-j_\sigma}l- m^*(l)n_\sigma^{-1}\right| \leq \frac{1}{2n_\sigma}$, for any $k\in I_\sigma$, which follows immediately from the definition of $m^*(l)$, to obtain
	\begin{align*}
		&\left|A\psi_{j_\sigma, l}^{(per)}(x) - Aw_{j_\sigma, m^*(l)}^{(per)}(x)\right|\\
		&= \left|\int_{-\infty}^\infty h(u-x)2^{j_\sigma/2}\psi(2^{j_\sigma}(u-2^{-j_\sigma}l))du - \int_{-\infty}^\infty h(u-x)2^{j_\sigma/2}\psi(2^{j_\sigma}(u-m^*(l)n_\sigma^{-1}))du\right|\\
		&=\left|\int_{-\infty}^\infty \left[h(2^{-j_\sigma}v+2^{-j_\sigma}l-x)-h(2^{-j_\sigma}v+m^*(l)n_\sigma^{-1}-x)\right]2^{-j_\sigma/2}\psi(v)dv\right|\\
		&= 2^{-3j_\sigma/2}\left|\int_{-\infty}^\infty \left[h'(2^{-j_\sigma}v+2^{-j_\sigma}l-x)-h'(2^{-j_\sigma}v+m^*(l)n_\sigma^{-1}-x)\right]\theta(v)dv\right|\\
		&\leq 2^{-3j_\sigma/2} C''\left|2^{-j_\sigma}l- m^*(l)n_\sigma^{-1}\right| \int_{-\infty}^\infty|\theta(v)|dv \leq  C''' \frac{2^{-3j_\sigma/2}}{2n_\sigma}.
	\end{align*}
	It follows that
	\[
	\inf_{0\leq l\leq 2^{j_\sigma}-1} \frac{\Vert A\psi_{j_\sigma, l}^{(per)} - Aw_{j_\sigma, m^*(l)}^{(per)}\Vert_\Y^2}{\Vert A\psi_{j_\sigma,0}^{(per)}\Vert_\Y^2} \leq C'''' \frac{n_\sigma^{-2}}{2^{3j_\sigma}\Vert A\psi_{j_\sigma,0}^{(per)}\Vert_\Y^2} \to 0,
	\]
	since $n_\sigma\to \infty$ and the denominator converges to a positive constant.

\end{proof}

\subsection{Proofs for section \ref{subsec:results_case2}}
Techniques used in the following proofs are inspired by \cite{laurentloubes2012}. Note that here we only consider the case that $\X$ and $\Y$ are complex spaces. The proofs for the case that they are real is analogous.

\subsubsection{Proof of the nonasymptotic upper bound}

\begin{proof}[Proof of Theorem \ref{thm:case2_upperbound_nonasymp}]
	Define the test
	\begin{equation}\label{eq:chisquaredtest}
		\phi_{\alpha,\sigma}(y_\sigma) =  \mathbbm{1}\left\{T_\sigma > t_{1-\alpha,\sigma} \right\},
	\end{equation}
	where $T_\sigma:=T_\sigma(Y_\sigma):= \sum_{k\in I_\sigma} |\langle Y_\sigma, v_k\rangle_\Y|^2$, and $t_{1-\alpha,\sigma}$ is the $(1-\alpha)$-quantile of $T_\sigma$ (which follows a generalized $\chi^2$-distribution) under $H_0$. Thus, by its very definition, $\phi_{\alpha,\sigma}$ is a level $\alpha$ test. We need to show that, if $\nu_\sigma$ is large enough, for any $f\in \mathcal{F}_\sigma^L(\nu_\sigma)$ 
	\begin{equation}\label{eq:gen_chisquared_powerbound}
		\Prob_f\left(T_\sigma \leq  t_{1-\alpha,\sigma} \right) \leq \delta-\alpha.
	\end{equation}
	We aim to show that asymptotically $t_{1-\alpha,\sigma}\leq t_{\delta-\alpha,\sigma}(f)$ whenever $\nu_\sigma\geq \sqrt{2} d_\alpha(\delta)\sigma\sqrt{\Vert \Xi_{I_\sigma}\Vert_F}$, where $t_{\delta-\alpha,\sigma}(f)$ denotes the $\delta-\alpha$ quantile of $T_\sigma$ when $f$ is the true underlying signal. First, we need to discuss the distribution of $T_\sigma$.
	
	For $f\in \mathcal{F}_\sigma^L(\nu_\sigma)$, the random vector $(\langle Y_\sigma, v_k\rangle_\Y)_{k\in I_\sigma}$ is normally distributed with with mean vector $m_\sigma= (\lambda_k\langle f,u_k\rangle_\X)_{k\in I_\sigma}$ and covariance matrix $2\sigma^2\Xi_{I_\sigma}$. Since $\Xi_\sigma$ is Hermitian and positive definite by assumption, it can be decomposed as 
	\[
	\Xi_{I_\sigma} = U_\sigma S_\sigma U_\sigma^H,
	\]
	where $U$ is unitary and $S_\sigma$ is a diagonal matrix containing the (real and positive) eigenvalues $(s_k)_{k\in I_\sigma}$ of $\Xi_\sigma$. It follows that  the random vector $(\langle Y_\sigma, v_k\rangle_\Y)_{k\in I_\sigma}$ can be written as
	$\sqrt{2}\sigma U_\sigma S_\sigma^{1/2} Z_\sigma + m_\sigma$ for some $Z_\sigma\sim \mathcal{CN}(0, \id_{N_\sigma})$ and thus,
	
	\begin{align*}
		T_\sigma &= (\sqrt{2}\sigma U_\sigma S_\sigma^{1/2} Z_\sigma + m_\sigma)^H(\sqrt{2}\sigma U_\sigma S_\sigma^{1/2} Z_\sigma + m_\sigma)\\
		&=2\sigma^2(Z_\sigma+ (\sqrt{2}\sigma)^{-1} S_\sigma^{-1/2} U_\sigma^H m_\sigma)^H S_\sigma(Z_\sigma+ (\sqrt{2}\sigma)^{-1} S_\sigma^{-1/2} U_\sigma^H m_\sigma)\\
		&= 2\sigma^2 \sum_{k\in I_\sigma} s_k \left|Z_{\sigma,k} -  \frac{1}{\sqrt{2}}\tilde{m}_{\sigma,k} \right|^2\\
		&=2 \sigma^2\sum_{k\in I_\sigma}  s_k \left[\left(\Re(Z_{\sigma,k}) -  \frac{1}{\sqrt{2}}\Re(\tilde{m}_{\sigma,k})\right)^2 + \left(\Im(Z_{\sigma,k}) -  \frac{1}{\sqrt{2}}\Im(\tilde{m}_{\sigma,k})\right)^2\right]\\
		&=\sigma^2  \sum_{k\in I_\sigma} s_k \left[(Z_k' - \Re(\tilde{m}_{\sigma,k}))^2 + (Z_k'' - \Im(\tilde{m}_{\sigma,k}))^2\right],
	\end{align*}
	where $\tilde{m}_\sigma=\sigma^{-1} S_\sigma^{-1/2} U_\sigma^H m_\sigma$ and  $Z',Z'' \iid \mathcal{N}(0,\id_{N_\sigma})$. In other words, $T_\sigma$ is the sum of $2N_\sigma$ weighted non-central chi-squared random variables. Note that
	\[
	\sigma^2\sum_{k\in I_\sigma} s_k |\tilde{m}_{\sigma,k}|^2 = (S_\sigma^{-1/2} U_\sigma^H m_\sigma)^H S (S_\sigma^{-1/2} U_\sigma^H m_\sigma) = m_\sigma^H m_\sigma = \sum_{k\in I_\sigma} |\lambda_k\langle f,u_k\rangle_\X|^2.
	\]
	
	\paragraph{Upper bound for $t_{1-\alpha,\sigma}$.} Under $H_0$ we have that $T_\sigma= \sigma^2\sum_{k\in I_\sigma} s_k \left[(Z_k')^2 + (Z_k'')^2\right]$. It follows from Lemma 1 from \cite{laurentmassart2000} that for any $t>0$
	\[
	\Prob_0\left(T_\sigma > 2\sigma^2\sum_{k\in I_\sigma} s_k + 2\sigma^2\sqrt{2t\sum_{k\in I_\sigma} s_k^2} + \sigma^2t\left(\sup_{k\in I_\sigma} s_k\right) \right)\leq \exp(-t),
	\]
	and thus, setting $t=\log(1/\alpha)$, we have
	\[
	t_{1-\alpha,\sigma} \leq  2\sigma^2\sum_{k\in I_\sigma} s_k + 2\sqrt{2\log(1/\alpha)}\sigma^2\Vert \Xi_{I_\sigma}\Vert_F + \log(1/\alpha)\sigma^2\sup_{k\in I_\sigma} s_k .
	\]
	
	\paragraph{Lower bound for $t_{\delta-\alpha,\sigma}(f)$.} We use Lemma 2 from \cite{laurentloubes2012}, which yields that for any $t>0$ 
	\[
	\Prob_f\left(T_\sigma\leq \E T_\sigma - 2\sqrt{2t}\sqrt{\sigma^4\sum_{k\in I_\sigma} s_k^2 + \sigma^4\sum_{k\in I_\sigma} s_k^2  |\tilde{m}_{\sigma,k}|^2}\right)\leq \exp(-t).
	\]
	Setting $t=\log(1/(\delta-\alpha))$, it follows that
	\begin{multline*}
		t_{\delta-\alpha,\sigma}(f) \geq \E T_\sigma -2\sqrt{2\log(1/(\delta-\alpha))}\sqrt{\sigma^4\Vert \Xi_{I_\sigma}\Vert_F^2 +\sigma^4 \sum_{k\in I_\sigma} s_k^2  |\tilde{m}_{\sigma,k}|^2}\\
		\geq 2\sigma^2\sum_{k\in I_\sigma} s_k + \sum_{k\in I_\sigma} |\lambda_k\langle f,u_k\rangle_\X|^2-2\sqrt{2\log(1/(\delta-\alpha))}\left[\sigma^2\Vert \Xi_{I_\sigma}\Vert_F + \sigma \sqrt{\Vert \Xi_{I_\sigma}\Vert_F}\sqrt{\sum_{k\in I_\sigma} |\lambda_k\langle f,u_k\rangle_\X|^2}\right],
	\end{multline*}
	where we used that $\sup_{k\in I_\sigma} s_k \leq \sqrt{\sum_{k\in I_\sigma} s_k^2}  = \Vert \Xi_{I_\sigma}\Vert_F$.
	\paragraph{Comparing the bounds.} It follows that $t_{1-\alpha,\sigma}\leq  t_{\delta-\alpha,\sigma}(f)$ is true when
	\begin{multline*}
		\sum_{k\in I_\sigma} |\lambda_k\langle f,u_k\rangle_\X|^2 - 2\sigma\sqrt{2\log(1/(\delta-\alpha))}\sqrt{\Vert \Xi_{I_\sigma}\Vert_F} \sqrt{\sum_{k\in I_\sigma} |\lambda_k\langle f,u_k\rangle_\X|^2}\\
		- \left(2\sqrt{2\log(1/(\delta-\alpha))}+2\sqrt{2\log(1/\alpha)}\right) \sigma^2 \Vert \Xi_{I_\sigma}\Vert_F -  \log(1/\alpha)\sigma^2 \sup_{k\in I_\sigma} s_k \geq 0,
	\end{multline*}
	which holds when
	\[
	\sqrt{\sum_{k\in I_\sigma} |\lambda_k\langle f,u_k\rangle_\X|^2} \geq \sqrt{2}\left[\sqrt{\log\frac{1}{\delta-\alpha}} +\left(\log\frac{1}{\alpha(\delta-\alpha)}+\sqrt{2\log\frac{1}{\delta-\alpha}}+\sqrt{2\log\frac{1}{\alpha}}\right)^{1/2}\right] \sigma\sqrt{\Vert\Xi_{I_\sigma}\Vert_F}.
	\]
\end{proof}

\subsubsection{Proof of the non-asymptotic lower bound}

\begin{proof}[Proof of Theorem \ref{thm:case2_lowerbound_nonasymp}] 
	The matrix $\tilde{\Xi}_{I_\sigma}$ given by $(\tilde{\Xi}_{I_\sigma})_{k,k'}= \langle \tilde{v}_k,\tilde{v}_{k'}\rangle_\Y$ is  Hermitian and positive definite, and thus, and we have the decompositions 
	\[
	\tilde{\Xi}_{I_\sigma} = \tilde{U}_\sigma \tilde{S}_\sigma \tilde{U}_\sigma^H, \quad \tilde{\Xi}_{I_\sigma}^{-1} = \tilde{U}_\sigma \tilde{S}_\sigma^{-1} \tilde{U}_\sigma^H,
	\]
	where $\tilde{U}_\sigma$ is unitary and $\tilde{S}_\sigma$ is a diagonal matrix with real and positive entries $(\tilde{s}_k)_{k\in I_\sigma}$ on its diagonal. The proof of the lower bound has the same core idea as the proof of Theorem \ref{thm:lowerbound_case1}: We start by defining a prior distribution on the set $\mathcal{F}_\sigma^L(\nu_\sigma)$. Let $w=(w_k)_{k\in I_\sigma}$ be a vector with $w_k\in \{-1,1\}$ for all $k$, and define
	\[
	\tilde{w} = \frac{\nu_\sigma}{\Vert \tilde{\Xi}_{I_\sigma}^{-1}\Vert_F} \tilde{U}_\sigma \tilde{S}_\sigma^{-1} w,
	\]
	and
	\[
	f_w= \sum_{k\in I_\sigma} \overline{\tilde{w}_k}\lambda_k^{-1} u_k.
	\]
	Note that 
	\[
	\sum_{k\in I_\sigma} |\lambda_k\langle f_w, u_k\rangle_\X|^2=\nu_\sigma^2,
	\]
	and thus, indeed $f_w\in \mathcal{F}_\sigma^L(\nu_\sigma)$. As in the proof of Theorem \ref{thm:lowerbound_case1} we get the likelihood ratio
	\[
	\frac{d\Prob_{f_w}}{d\Prob_0}(\tilde{y}_\sigma)= \exp\left(\frac{1}{\sigma}\sum_{j\in \N} \left[\Re(\langle Af_w, e_j\rangle_\Y) \tilde{y}_{\sigma,2j-1} +\Im(\langle Af_w, e_j\rangle_\Y) \tilde{y}_{\sigma,2j}\right] -\frac{\Vert Af_w\Vert_\Y^2}{2\sigma^2} \right).
	\]
	Note that $Af_w= \sum_{k\in I_\sigma} \overline{\tilde{w}_k} \tilde{v}_k$, an thus,
	\[
	\Vert Af_w\Vert_\Y^2 = \sum_{k,l\in I_\sigma} \overline{\tilde{w}_k} \tilde{w}_l \langle\tilde{v}_k, \tilde{v}_l\rangle_\Y = \frac{\nu_\sigma^2}{\Vert \tilde{\Xi}_\sigma^{-1}\Vert_F^2} w^T \tilde{S}_\sigma^{-1} \tilde{U}_\sigma^H \tilde{\Xi}_{I_\sigma} \tilde{U}_\sigma\tilde{S}_\sigma^{-1}  w = \frac{\nu_\sigma^2}{\Vert \tilde{\Xi}_{I_\sigma}^{-1}\Vert_F^2}\sum_{k\in I_\sigma} \frac{1}{\tilde{s}_k}.
	\]
	Let $\tilde{v}^{(j)}$ be the vector with entries $\tilde{v}^{(j)}_k = \langle \tilde{v}_k, e_j\rangle_\Y$ for $k\in I_\sigma$. Then
	\[
	\langle Af_w, e_j\rangle_\Y = \sum_{k\in I_\sigma}\overline{\tilde{w}_k}\langle \tilde{v}_k, e_j\rangle_\Y = \frac{\nu_\sigma}{\Vert \tilde{\Xi}_{I_\sigma}^{-1}\Vert_F} w^T \tilde{S}_\sigma^{-1}\tilde{U}_\sigma^H \tilde{v}^{(j)}= \frac{\nu_\sigma}{\Vert \tilde{\Xi}_{I_\sigma}^{-1}\Vert_F}\sum_{k\in I_\sigma}w_k \left(\tilde{S}_\sigma^{-1}\tilde{U}_\sigma^H \tilde{v}^{(j)}\right)_k,
	\]
	and it follows that
	\[
	\sum_{j\in \N} \left[\Re(\langle Af_w, e_j\rangle_\Y) \tilde{y}_{\sigma,2j-1} +\Im(\langle Af_w, e_j\rangle_\Y) \tilde{y}_{\sigma,2j}\right] = \frac{\nu_\sigma}{\Vert \tilde{\Xi}_{I_\sigma}^{-1}\Vert_F}\sum_{k\in I_\sigma}w_k Z_{\sigma,k},
	\]
	where $Z_\sigma = (Z_{\sigma,k})_{k\in I_\sigma}$, with $Z_{\sigma,k}= \sum_{j\in \N}\left[\Re \left(\tilde{S}_\sigma^{-1}\tilde{U}_\sigma^H \tilde{v}^{(j)}\right)_k\tilde{y}_{\sigma,2j-1} + \Im\left(\tilde{S}_\sigma^{-1}\tilde{U}_\sigma^H \tilde{v}^{(j)}\right)_k\tilde{y}_{\sigma,2j}\right]$, is a normally distributed random vector with mean $0$ and covariance matrix $\Sigma$ given by
	\begin{align*}
		\Sigma_{k,l} &= \sum_{j\in \N} \Re\left[\overline{\left(\tilde{S}_\sigma^{-1}\tilde{U}_\sigma^H \tilde{v}^{(j)}\right)_k}\left(\tilde{S}_\sigma^{-1}\tilde{U}_\sigma^H \tilde{v}^{(j)}\right)_l\right]\\
		&= \sum_{j\in \N} \Re\left[\overline{\left(\sum_{m\in I_\sigma}\left(\tilde{S}_\sigma^{-1}\tilde{U}_\sigma^H\right)_{k,m} \langle \tilde{v}_m, e_j\rangle_\Y\right)}\left(\sum_{m'\in I_\sigma}\left(\tilde{S}_\sigma^{-1}\tilde{U}_\sigma^H\right)_{l,m'} \langle \tilde{v}_{m'}, e_j\rangle_\Y\right)\right]\\
		&= \Re\left[\sum_{m,m'\in \N} \overline{\left(\tilde{S}_\sigma^{-1}\tilde{U}_\sigma^H\right)_{k,m}}\left(\tilde{S}_\sigma^{-1}\tilde{U}_\sigma^H\right)_{l,m'} \langle \tilde{v}_{m'}, \tilde{v}_m\rangle_\Y\right]=\Re\left(\tilde{S}_\sigma^{-1} \tilde{U}_\sigma^H \tilde{\Xi}_\sigma \tilde{U}_\sigma\tilde{S}_\sigma^{-1}\right)_{k,l} = \frac{1}{\tilde{s}_k}\delta_{k,l}.
	\end{align*}
	In other words, the random variables $Z_{\sigma,k}$ are independent and $Z_{\sigma,k}\sim\mathcal{N}(0,1/\tilde{s}_k)$ for $k\in I_\sigma$. Thus, under $H_0$, we have
	\[
	\frac{d\Prob_{f_w}}{d\Prob_0}(\tilde{y}_\sigma) =  \exp\left[-\frac{\nu_\sigma^2}{2\sigma^2\Vert \tilde{\Xi}_{I_\sigma}^{-1}\Vert_F^2}\sum_{k\in I_\sigma} \frac{1}{\tilde{s}_k}\right]\prod_{k\in I_\sigma}\exp\left[\frac{\nu_\sigma}{\sigma\Vert \tilde{\Xi}_{I_\sigma}^{-1}\Vert_F}\cdot \frac{w_k}{\sqrt{\tilde{s}_k}} Z_{\sigma,k}'\right],
	\]
	where $Z_{\sigma,k}'\iid \mathcal{N}(0,1)$ for $k\in I_\sigma$.
	
	Now, assume that $\hat{w}_k$, $k\in I_\sigma$ are independent Rademacher variables (which means that $\Prob(\hat{w}_k=1)=\Prob(\hat{w}_k=-1)=1/2$ for any $k$), that are also independent from $\tilde{y}_\sigma$, and let $\hat{w}=(\hat{w}_k)_{k\in I_\sigma}$ be the corresponding random vector. We denote by $\pi_\sigma$ the (finitely supported) distribution of the random function $f_{\hat{w}}$  on $\mathcal{F}_\sigma^L(\nu_\sigma)$. As in the proof of Theorem \ref{thm:lowerbound_case1} we have
	\[
	\gamma_\sigma= \inf_{\phi \in \Phi_\sigma}\left[\alpha_\sigma(\phi)+ \sup_{f\in\mathcal{F}_\sigma^L(\nu_\sigma)}\beta_\sigma(\phi, f)\right] = 1-  \frac{1}{2}\E_0 \left|\E_{\pi_\sigma}\frac{d\Prob_f}{d\Prob_0}(Y_\sigma)-1\right|,
	\]
	Note that 
	\[
	\E_0 \E_{\pi_\sigma}\frac{d\Prob_f}{d\Prob_0}(\tilde{y}_\sigma) =  \E_{\pi_\sigma} \E_0\frac{d\Prob_f}{d\Prob_0}(\tilde{y}_\sigma) =1,
	\]
	and it follows that
	\[
	\left(\E_0 \left|\E_{\pi_\sigma}\frac{d\Prob_f}{d\Prob_0}(\tilde{y}_\sigma)-1\right|\right)^2\leq \E_0\left(\E_{\pi_\sigma}\frac{d\Prob_f}{d\Prob_0}(\tilde{y}_\sigma)\right)^2-1,
	\]
	and thus, 
	\[
	\gamma_\sigma(\nu_\sigma) \geq 1-\frac{1}{2}\left(\E_0\left(\E_{\pi_\sigma}\frac{d\Prob_f}{d\Prob_0}(\tilde{y}_\sigma)\right)^2-1\right)^{1/2}.
	\]
	We have
	\[
	\E_{\pi_\sigma}\frac{d\Prob_f}{d\Prob_0}(\tilde{y}_\sigma) = \exp\left[-\frac{\nu_\sigma^2}{2\sigma^2\Vert \tilde{\Xi}_{I_\sigma}^{-1}\Vert_F^2}\sum_{k\in I_\sigma} \frac{1}{\tilde{s}_k}\right]\prod_{k\in I_\sigma}\cosh\left[\frac{\nu_\sigma}{\sigma\sqrt{\tilde{s}_k}\Vert \tilde{\Xi}_{I_\sigma}^{-1}\Vert_F} Z_{\sigma,k}'\right].
	\]
	It follows that
	\[
	\E_0\left(\E_{\pi_\sigma}\frac{d\Prob_f}{d\Prob_0}(\tilde{y}_\sigma)\right)^2 = \prod_{k\in I_\sigma} \cosh\left[\frac{\nu_\sigma^2}{\sigma^2\tilde{s}_k \Vert \tilde{\Xi}_{I_\sigma}^{-1}\Vert_F^2} \right]\leq  \prod_{k\in I_\sigma} \exp\left[\frac{\nu_\sigma^4}{\sigma^4\tilde{s}_k \Vert \tilde{\Xi}_{I_\sigma}^{-1}\Vert_F^4}\right] = \exp\left[\frac{\nu_\sigma^4}{\sigma^4 \Vert \tilde{\Xi}_{I_\sigma}^{-1}\Vert_F^2}\right],
	\]
	where we used that $\E\cosh^2(tX) = \exp(t^2)\cosh(t^2)$ for any $t\in \R$ and $X\sim\mathcal{N}(0,1)$ and that $\cosh(t)\leq \exp(t^2/2)$ for any $t\in \R$. The claim follows immediately.
	
\end{proof}

\subsubsection{Remaining proofs}

\begin{proof}[Proof of Corollary \ref{cor:case2_asymp}] The second part follows immediately from Theorem \ref{thm:case2_lowerbound_nonasymp}. For the first part, note that $\gamma_\sigma(\nu_\sigma)\leq \gamma_{\sigma,\alpha}(\nu_\sigma)$ for any $\alpha$. Thus, the first part follows immediately from Theorem \ref{thm:case2_upperbound_nonasymp}.
\end{proof}

\begin{proof}[Proof of Lemma \ref{lm:rate_for_constant_norm}]
	\textit{(1)} Let $z= (z_k)_{k\in I_\sigma}$ be a non-zero complex vector. Then it follows from the fact that $(v_k)_{k\in I_\sigma}$ is a Riesz sequence that
	\[
	z^H\Xi_{I_\sigma}z = \left\Vert \sum_{k\in I_\sigma}\overline{z_k}v_k\right\Vert_\Y\geq \left(C\sum_{k\in I_\sigma}|\overline{z_k}|^2\right)^{1/2}>0,
	\]
	for some constant $C>0$. The proof for $\tilde{\Xi}_{I_\sigma}$ is analogous.\\
	\textit{(2)} The results of Theorem \ref{thm:case2_upperbound_nonasymp} and \ref{thm:case2_lowerbound_nonasymp} (which can be applied since $\Xi_{I_\sigma}$ and $\tilde{\Xi}_{I_\sigma}$ are positive definite) imply that $\Vert \tilde{\Xi}_\sigma^{-1}\Vert_F = O(\Vert \Xi_{I_\sigma}\Vert_F)$. It remains to show that $\Vert \Xi_{I_\sigma}\Vert_F \leq C \Vert \tilde{\Xi}_\sigma^{-1}\Vert_F$ for some constant $C>0$.  We have
	\[
	\Vert \Xi_{I_\sigma}\Vert_F = \Vert \Xi_{I_\sigma} \tilde{\Xi}_{I_\sigma} \tilde{\Xi}_{I_\sigma}^{-1} \Vert_F \leq \Vert \Xi_{I_\sigma} \tilde{\Xi}_{I_\sigma}\Vert_2\Vert \tilde{\Xi}_{I_\sigma}^{-1}\Vert_F,
	\]
	where $\Vert \Xi_{I_\sigma} \tilde{\Xi}_{I_\sigma}\Vert_2 = \max_{\Vert z\Vert_2=1} \Vert \Xi_{I_\sigma}\tilde{\Xi}_{I_\sigma} z \Vert_2$, where $\Vert \cdot \Vert_2$ denotes the euclidean norm on $\C^{N_\sigma}$. Now let $z= (z_k)_{k\in I_\sigma}$ be a complex vector with $\Vert z\Vert_2=1$. Recall that, since $(v_k)_{k\in I}$ and $(\tilde{v}_k)_{k\in I}$ are Riesz sequences, they are also frames of their respective spans. It follows that
	\begin{align*}
		\Vert \Xi_{I_\sigma}\tilde{\Xi}_{I_\sigma} z \Vert_2^2 &= \sum_{k\in I_\sigma} \left|\sum_{l\in I_\sigma}\left(\sum_{j\in I_\sigma} \langle v_k,v_j\rangle_\Y  \langle \tilde{v}_j,\tilde{v}_{l}\rangle_\Y\right) z_{l} \right|^2\leq \sum_{k\in I} \left|\left\langle\sum_{j\in I_\sigma} \left\langle \sum_{l\in I_\sigma} \overline{z_l}\tilde{v}_{l}, \tilde{v}_j\right\rangle_\Y v_j,  v_k\right\rangle_\Y    \right|^2\\
		&\leq C \left\Vert \sum_{j\in I_\sigma} \left\langle \sum_{l\in I_\sigma} \overline{z_l}\tilde{v}_{l}, \tilde{v}_j\right\rangle_\Y v_j\right\Vert_\Y^2\leq C' \sum_{j\in I_\sigma}\left|\left\langle \sum_{l\in I_\sigma} \overline{z_l}\tilde{v}_{l}, \tilde{v}_j\right\rangle_\Y\right|^2 \\
		&\leq C' \sum_{j\in I}\left|\left\langle \sum_{l\in I_\sigma} \overline{z_l}\tilde{v}_{l}, \tilde{v}_j\right\rangle_\Y\right|^2\leq C'' \left\Vert \sum_{l\in I_\sigma} \overline{z_l}\tilde{v}_{l}\right\Vert_\Y^2 \leq C''' \sum_{l\in I_\sigma}|\overline{z_l}|^2 = C''',
	\end{align*}
	which concludes this proof.\\
	\textit{(3)} Note that there are constants $C_1, C_2>0$, such that $C_1\leq \Vert v_k\Vert_\Y\leq C_2$, for any $k\in I$, since $(v_k)_{k\in I}$ is a Riesz sequence. It follows that
	\[
	\Vert \Xi_{I_\sigma}\Vert_F^2= \sum_{k\in I_\sigma}\sum_{k'\in I_\sigma} |\langle v_k, v_{k'}\rangle_\Y|^2 \geq \sum_{k\in I_\sigma}|\langle v_k, v_{k}\rangle_\Y|^2 = \sum_{k\in I_\sigma} \Vert v_k\Vert_\Y^4 \geq C_1^4N_\sigma,
	\]
	and
	\[
	\Vert \Xi_{I_\sigma}\Vert_F^2\leq \sum_{k\in I_\sigma}\sum_{k'\in I} |\langle v_k, v_{k'}\rangle_\Y|^2 \leq C  \sum_{k\in I_\sigma}\Vert v_k\Vert_\Y^2 \leq CC_2^2N_\sigma.
	\]
	The claim follows.
\end{proof}
\bibliographystyle{plain}
\bibliography{ref}

\begin{thebibliography}{10}

\bibitem{abramovich1998vwd}
F.~Abramovich and B.~W. Silverman.
\newblock Wavelet decomposition approaches to statistical inverse problems.
\newblock {\em Biometrika}, 85(1):115--129, 1998.

\bibitem{ab06}
Anestis Antoniadis and J\'{e}remie Bigot.
\newblock Poisson inverse problems.
\newblock {\em Ann. Statist.}, 34(5):2132--2158, 2006.

\bibitem{autin2019maxiset}
F.~Autin, M.~Clausel, J.-M. Freyermuth, and C.~Marteau.
\newblock Maxiset point of view for signal detection in inverse problems.
\newblock {\em Math. Methods Statist.}, 28(3):228--242, 2019.

\bibitem{bhmr07}
N.~Bissantz, T.~Hohage, A.~Munk, and F.~Ruymgaart.
\newblock Convergence rates of general regularization methods for statistical
  inverse problems and applications.
\newblock {\em SIAM J. Numer. Anal.}, 45(6):2610--2636, 2007.

\bibitem{bmp09}
Cristina Butucea, Catherine Matias, and Christophe Pouet.
\newblock Adaptive goodness-of-fit testing from indirect observations.
\newblock {\em Ann. Inst. Henri Poincar\'{e} Probab. Stat.}, 45(2):352--372,
  2009.

\bibitem{daubechies1992ten}
Ingrid Daubechies.
\newblock {\em Ten lectures on wavelets}, volume~61 of {\em CBMS-NSF Regional
  Conference Series in Applied Mathematics}.
\newblock Society for Industrial and Applied Mathematics (SIAM), Philadelphia,
  PA, 1992.

\bibitem{am20}
Miguel del \'{A}lamo and Axel Munk.
\newblock Total variation multiscale estimators for linear inverse problems.
\newblock {\em Inf. Inference}, 9(4):961--986, 2020.

\bibitem{donoho1995wvd}
David~L. Donoho.
\newblock Nonlinear solution of linear inverse problems by wavelet-vaguelette
  decomposition.
\newblock {\em Appl. Comput. Harmon. Anal.}, 2(2):101--126, 1995.

\bibitem{ebner2020}
Andrea Ebner, J\"{u}rgen Frikel, Dirk Lorenz, Johannes Schwab, and Markus
  Haltmeier.
\newblock Regularization of inverse problems by filtered diagonal frame
  decomposition.
\newblock {\em Appl. Comput. Harmon. Anal.}, 62:66--83, 2023.

\bibitem{ehn96}
Heinz~W. Engl, Martin Hanke, and Andreas Neubauer.
\newblock {\em Regularization of inverse problems}, volume 375 of {\em
  Mathematics and its Applications}.
\newblock Kluwer Academic Publishers Group, Dordrecht, 1996.

\bibitem{enikeeva2020}
Farida Enikeeva, Axel Munk, Markus Pohlmann, and Frank Werner.
\newblock Bump detection in the presence of dependency: does it ease or does it
  load?
\newblock {\em Bernoulli}, 26(4):3280--3310, 2020.

\bibitem{enikeeva2018}
Farida Enikeeva, Axel Munk, and Frank Werner.
\newblock Bump detection in heterogeneous {G}aussian regression.
\newblock {\em Bernoulli}, 24(2):1266--1306, 2018.

\bibitem{gn20}
Matteo Giordano and Richard Nickl.
\newblock Consistency of {B}ayesian inference with {G}aussian process priors in
  an elliptic inverse problem.
\newblock {\em Inverse Problems}, 36(8):085001, 35, 2020.

\bibitem{hlm22}
Markus Haltmeier, Housen Li, and Axel Munk.
\newblock A variational view on statistical multiscale estimation.
\newblock {\em Annu. Rev. Stat. Appl.}, 9:343--372, 2022.

\bibitem{hanke2017}
Martin Hanke.
\newblock {\em A taste of inverse problems}.
\newblock Society for Industrial and Applied Mathematics (SIAM), Philadelphia,
  PA, 2017.
\newblock Basic theory and examples.

\bibitem{hbm05}
Hajo Holzmann, Nicolai Bissantz, and Axel Munk.
\newblock Density testing in a contaminated sample.
\newblock {\em J. Multivariate Anal.}, 98(1):57--75, 2007.

\bibitem{hubmer2021}
Simon Hubmer and Ronny Ramlau.
\newblock Frame decompositions of bounded linear operators in {H}ilbert spaces
  with applications in tomography.
\newblock {\em Inverse Problems}, 37(5):Paper No. 055001, 30, 2021.

\bibitem{ingster2014}
Yu. Ingster, B.~Laurent, and C.~Marteau.
\newblock Signal detection for inverse problems in a multidimensional
  framework.
\newblock {\em Math. Methods Statist.}, 23(4):279--305, 2014.

\bibitem{ingster1993}
Yu.~I. Ingster.
\newblock Asymptotically minimax hypothesis testing for nonparametric
  alternatives. {I-III}.
\newblock {\em Math. Methods Statist.}, 2(4):249--268, 1993.

\bibitem{ingster2003}
Yu.~I. Ingster and I.~A. Suslina.
\newblock {\em Nonparametric goodness-of-fit testing under {G}aussian models},
  volume 169 of {\em Lecture Notes in Statistics}.
\newblock Springer-Verlag, New York, 2003.

\bibitem{ingster2012}
Yuri~I. Ingster, Theofanis Sapatinas, and Irina~A. Suslina.
\newblock Minimax signal detection in ill-posed inverse problems.
\newblock {\em Ann. Statist.}, 40(3):1524--1549, 2012.

\bibitem{j99}
Iain~M. Johnstone.
\newblock Wavelet shrinkage for correlated data and inverse problems:
  adaptivity results.
\newblock {\em Statist. Sinica}, 9(1):51--83, 1999.

\bibitem{johnstone2004}
Iain~M. Johnstone, G\'{e}rard Kerkyacharian, Dominique Picard, and Marc
  Raimondo.
\newblock Wavelet deconvolution in a periodic setting.
\newblock {\em J. R. Stat. Soc. Ser. B Stat. Methodol.}, 66(3):547--573, 2004.

\bibitem{js97}
Iain~M. Johnstone and Bernard~W. Silverman.
\newblock Wavelet threshold estimators for data with correlated noise.
\newblock {\em J. Roy. Statist. Soc. Ser. B}, 59(2):319--351, 1997.

\bibitem{klosowski2017}
Jakob Klosowski and Jens Frahm.
\newblock Image denoising for real-time mri.
\newblock {\em Magnetic resonance in medicine}, 77(3):1340--1352, 2017.

\bibitem{kretschmann2022}
Remo Kretschmann, Daniel Wachsmuth, and Frank Werner.
\newblock Optimal regularized hypothesis testing in statistical inverse
  problems.
\newblock arXiv: 2212.12897, 2022.

\bibitem{laurent2011directindirect}
B.~Laurent, J.-M. Loubes, and C.~Marteau.
\newblock Testing inverse problems: a direct or an indirect problem?
\newblock {\em J. Statist. Plann. Inference}, 141(5):1849--1861, 2011.

\bibitem{laurentmassart2000}
B.~Laurent and P.~Massart.
\newblock Adaptive estimation of a quadratic functional by model selection.
\newblock {\em Ann. Statist.}, 28(5):1302--1338, 2000.

\bibitem{laurentloubes2012}
B\'{e}atrice Laurent, Jean-Michel Loubes, and Cl\'{e}ment Marteau.
\newblock Non asymptotic minimax rates of testing in signal detection with
  heterogeneous variances.
\newblock {\em Electron. J. Stat.}, 6:91--122, 2012.

\bibitem{l15}
Oleg Lepski.
\newblock Adaptive estimation over anisotropic functional classes via oracle
  approach.
\newblock {\em Ann. Statist.}, 43(3):1178--1242, 2015.

\bibitem{lifshits2012}
Mikhail Lifshits.
\newblock {\em Lectures on {G}aussian processes}.
\newblock SpringerBriefs in Mathematics. Springer, Heidelberg, 2012.

\bibitem{lm14}
J.~M. Loubes and C.~Marteau.
\newblock Goodness-of-fit testing strategies from indirect observations.
\newblock {\em J. Nonparametr. Stat.}, 26(1):85--99, 2014.

\bibitem{mallat}
St\'{e}phane Mallat.
\newblock {\em A wavelet tour of signal processing}.
\newblock Academic Press, Inc., San Diego, CA, 1998.

\bibitem{marteau2014general}
C.~Marteau and P.~Math\'{e}.
\newblock General regularization schemes for signal detection in inverse
  problems.
\newblock {\em Math. Methods Statist.}, 23(3):176--200, 2014.

\bibitem{m09}
Alexander Meister.
\newblock {\em Deconvolution problems in nonparametric statistics}, volume 193
  of {\em Lecture Notes in Statistics}.
\newblock Springer-Verlag, Berlin, 2009.

\bibitem{mnp21}
Fran\c{c}ois Monard, Richard Nickl, and Gabriel~P. Paternain.
\newblock Statistical guarantees for {B}ayesian uncertainty quantification in
  nonlinear inverse problems with {G}aussian process priors.
\newblock {\em Ann. Statist.}, 49(6):3255--3298, 2021.

\bibitem{munkstaudtwerner2020}
Axel Munk, Thomas Staudt, and Frank Werner.
\newblock Statistical foundations of nanoscale photonic imaging.
\newblock {\em Nanoscale Photonic Imaging}, pages 125--143, 2020.

\bibitem{natterer2001}
F.~Natterer.
\newblock {\em The mathematics of computerized tomography}.
\newblock B. G. Teubner, Stuttgart; John Wiley \& Sons, Ltd., Chichester, 1986.

\bibitem{ot21}
Taisuke Otsu and Luke Taylor.
\newblock Specification testing for errors-in-variables models.
\newblock {\em Econometric Theory}, 37(4):747--768, 2021.

\bibitem{proksch2018}
Katharina Proksch, Frank Werner, and Axel Munk.
\newblock Multiscale scanning in inverse problems.
\newblock {\em Ann. Statist.}, 46(6B):3569--3602, 2018.

\bibitem{rs16}
Kolyan Ray and Johannes Schmidt-Hieber.
\newblock Minimax theory for a class of nonlinear statistical inverse problems.
\newblock {\em Inverse Problems}, 32(6):065003, 29, 2016.

\bibitem{t00}
Alexandre Tsybakov.
\newblock On the best rate of adaptive estimation in some inverse problems.
\newblock {\em C. R. Acad. Sci. Paris S\'{e}r. I Math.}, 330(9):835--840, 2000.

\bibitem{vardi1985}
Y.~Vardi, L.~A. Shepp, and L.~Kaufman.
\newblock A statistical model for positron emission tomography.
\newblock {\em J. Amer. Statist. Assoc.}, 80(389):8--37, 1985.
\newblock With discussion.

\bibitem{wang2014}
Xinghui Wang and Shuhe Hu.
\newblock Weak laws of large numbers for arrays of dependent random variables.
\newblock {\em Stochastics}, 86(5):759--775, 2014.

\bibitem{wernerhohage2012}
Frank Werner and Thorsten Hohage.
\newblock Convergence rates in expectation for {T}ikhonov-type regularization
  of inverse problems with {P}oisson data.
\newblock {\em Inverse Problems}, 28(10):104004, 15, 2012.

\end{thebibliography}
\end{document}